\documentclass[a4paper,12pt,reqno]{amsart}
\usepackage{amssymb}
\usepackage{latexsym}
\usepackage{amsmath}
\usepackage[english]{babel}
\usepackage{amsfonts}
\usepackage{enumerate}
\usepackage[all]{xy}
\CompileMatrices

\usepackage[margin=3cm,footskip=2em,top=2.8cm]{geometry}

\usepackage{rotating}
\newlength{\widthuparrow}
\newcommand{\Coker}{\operatorname{Coker}}
\newcommand{\Hom}{\operatorname{Hom}}
\newcommand{\Ext}{\operatorname{Ext}}
\newcommand{\Tor}{\operatorname{Tor}}
\newcommand{\Ind}{\operatorname{Ind}}

\newcommand{\End}{\operatorname{End}}

\newcommand{\cha}{\operatorname{char}}
\newcommand{\ch}{\operatorname{ch}}
\newcommand{\di}{\operatorname{div}}
\newcommand{\Di}{\operatorname{Div}}

\newcommand{\sign}{\operatorname{sign}}
\usepackage{bbm} 
\newcommand\C{\mathbb C}
\newcommand\Z{\mathbb Z}
\newcommand\R{\mathbb R}

\newcommand\N{\mathbb N}
\newcommand\Q{\mathbb Q}

\newcommand\pf{\noindent {\bf Proof:  }}

\usepackage{amsthm}
\newtheorem{thm}{Theorem}[section]

\newtheorem{prop}[thm]{Proposition}
\newtheorem{cor}[thm]{Corollary}
\newtheorem{lemma}[thm]{Lemma}

\theoremstyle{definition}
\newtheorem{rem}[thm]{Remark}

\numberwithin{equation}{section}

\title{Sum formulas for reductive algebraic groups}
\author{Henning Haahr Andersen and Upendra Kulkarni}

\address{HHA: Department of Mathematics, University of Aarhus
, Building 530, Ny
Munkegade,
8000  Aarhus C, Denmark}
\email{mathha@imf.au.dk}

\address{UK: TIFR Centre, P.O. Box 1234, IISc Campus, Bangalore 560 012, India}
\email{upendra@math.tifrbng.res.in}

\begin{document}

\begin{abstract}
Let $V$ be a Weyl module either for a reductive algebraic group $G$ or
for the corresponding quantum group $U_q$. If $G$ is defined over 
a field of positive characteristic $p$, respectively if $q$ is a primitive $l$'th root of unity 
(in an arbitrary field) then $V$ has a Jantzen filtration $V=V^0
\supset V^1 \cdots \supset V^r = 0$. The sum of the positive terms in
this filtration satisfies a well known sum formula. 

If $T$ denotes a tilting module either for $G$ or $U_q$ then we can
similarly filter the space $\Hom_G(V,T)$, respectively 
$\Hom_{U_q}(V,T)$  and there is a sum
formula for the positive terms here as well.

We give an easy and unified proof of these two (equivalent) sum
formulas. Our approach is based on an Euler type identity which we
show holds without any restrictions on $p$ or $l$. In
particular, we get rid of previous such restrictions in the tilting
module case.
\end{abstract}

\maketitle

\section{Introduction}
Let $G$ be a reductive algebraic group over a field $k$ of prime
characteristic $p$.  The Weyl modules play a fundamental role in the
study of finite dimensional representations of $G$. One of the
important tools in investigating the structure of Weyl modules is
their Jantzen filtration. The sum of the characters of the filtration
terms obeys a sum formula analogous to the Verma module case,
\cite{MHG}. This formula was first proved by J. C. Jantzen some 30
years ago \cite{JCJ} with some mild restrictions on~$p$. Later the
first author \cite{Asum} gave another proof valid for all $p$ based on
the fact that Weyl modules are special cases of cohomology of line
bundles on the flag manifold for $G$ and exploring natural
homomorphisms between such cohomology modules.

More recently a similar formula \cite{Atiltsum} turned up in the
theory for tilting modules for~$G$. This time one filters the space of
homomorphisms from a Weyl module into a tilting module. However, the
proof in loc. cit. needs $p$ to be at least the Coxeter number for
$G$.

Both the above mentioned sum formulas are related to Ext-groups involving
integral versions of Weyl modules, see \cite{Kul1} and \cite{Atiltsum}. In this paper
we start out by proving an Euler-type formula for such Ext-groups using
techniques from \cite{Asum} and \cite{Kul1}. Then we are able to deduce the two 
sum formulas from this. In particular, our results work for all $p$. It also
reveals that the two cases are in fact equivalent.

Let $U_q$ denote the quantum group corresponding to $G$. When $q$ is a root
of unity (and $U_q$ is obtained via the Lusztig divided power construction)
there are completely analogous sum formulas for $U_q$. Our proof applies in
this case as well and it avoids the restrictions on the
order of $q$ in \cite{Atiltsum}.

We have taken the opportunity to recall the arguments from \cite{Asum}, \cite{Atiltsum} and 
\cite{Kul1} that we need. In this way our proof of the sum formulas for $G$ is completely 
self-contained relying only on basic facts on Weyl modules, cohomology on 
line bundles, and tilting modules (which can all be found in \cite{RAG}). In the
quantum case everything works in the same way and we have only given the statements 
in that case leaving the analogous proofs to the readers. 

Some of the results in this paper date back several years. At the meeting
AMS Scand 2000 in Odense, Denmark the second author gave a talk, 
``Ext groups and Jantzen's sum formula'' in which he presented the
Weyl module sum formula in terms of Ext-groups. This can be found in \cite{Kul1}, and it is 
also referred to in the preprint \cite{Kul2} where he proves the equivalence with the 
sum formula for tilting modules. Shortly after the appearance of this
preprint we realized how to give the
uniform proof presented below.

\section{Notation}
\subsection{Roots}
Throughout this paper $k$ will denote an algebraically closed field of 
characteristic $p > 0$ and $G$ will denote a reductive algebraic group over $k$. We 
choose a maximal torus $T$ in $G$ and a Borel subgroup $B$ containing $T$.
Then $R$ will be the root system for $(G,T)$. We fix a set of simple roots
$S$ in $R$ by requiring that the roots of $B$ are the corresponding negative
roots $-R^+$. The number of positive roots is called $N$. This is also the
dimension of the flag variety $G/B$.

The character group for $T$ (and $B$) is denoted $X$. We let $X^+$ be the
set of dominant characters, i.e., $X^+ = \{\lambda \in X \mid \langle \lambda,
\alpha^\vee \rangle\geq ~0 \text { for all } \alpha~\in~R^+\}$.

The Weyl group $W = N_G(T)/Z_G(T)$ for $G$ acts naturally on $X$. If $\alpha \in
R$ then the reflection $s_\alpha \in W$ corresponding to $\alpha$ is given by
$s_\alpha(\lambda) = \lambda - \langle\lambda, \alpha^\vee \rangle \alpha$ for
all $\lambda \in X$. We shall also use the `dot-action' defined by $w \cdot
\lambda = w(\lambda + \rho) - \rho, w \in W, \lambda \in X$. Here $\rho$ is the
half sum of the positive roots. 

Each element $w \in W$ is a product of simple 
reflections (reflections for simple roots) and we have the corresponding length 
function $l$ on $W$ taking $w$ into the minimal number of such simple reflections
needed to express $w$. The unique longest element in $W$ is denoted $w_0$. It has 
length $l(w_0) = N$.

\subsection{Weights}
If $M$ is a finite dimensional $T$-module and $\lambda \in X$ then the weight
space $M_\lambda$ is defined by $M_\lambda = \{m\in M\mid t m = \lambda(t)m
\text { for all } t\in T\}$. We say that $\lambda$ is a weight of $M$ if $M_\lambda \not =
0$. The character $\ch M$ is $\ch M =
\sum_{\lambda \in X} (\dim M_\lambda) e^\lambda \in \Z[X]$.

For each $\lambda \in X^+$ we have a Weyl module $\Delta (\lambda)$ for $G$ with highest weight 
$\lambda$. Its contragredient dual $\Delta (\lambda)^*$ is denoted 
$\nabla(-w_0\lambda)$. Note that then the dual Weyl module $\nabla(\mu)$ attached
to $\mu \in X^+$ has highest weight $\mu$ (because $w_0(\lambda)$ is the smallest weight
of $\Delta(\lambda)$).

\subsection{Cohomology modules}
Let $M$ be a finite dimensional $B$-module. Then we will write $H^0(M)$ for the
$G$-module $\Ind_B^G M$ induced by $M$. This is also the $0$-th cohomology (i.e.,
the set of global sections) for the vector bundle on $G/B$ associated with $M$.
More generally, we denote by $H^i(M)$ the $i$-th cohomology of this bundle, or
alternatively the value of the $i$-th right derived functor $R^i\Ind_B^G$ on $M$.
It is well known (as $G/B$ is a projective variety) that the cohomology $H^\bullet(M)$
is finite dimensional, and that $H^i(M) = 0 $ for $i>N$.

The Euler character of a $B$-module $M$ is given by
$$ \chi (M) = \sum_{\mu \in X} (-1)^i\ch (H^i(M)).$$
Note that $\chi $ is additive, i.e., if $0 \rightarrow M_1 \rightarrow M \rightarrow
M_2 \rightarrow 0$ is a short exact sequence of finite dimensional $B$-modules then 
$\chi (M) = \chi (M_1) + \chi (M_2)$.

In the following the cohomology modules $H^i(\lambda), \lambda \in X$ will
play a vital role. In particular, we recall that the Weyl modules above are special 
instances of such modules. Precisely, we have $\Delta (\lambda) \simeq H^N(w_0\cdot \lambda)$
for all $\lambda \in X^+$. Also $\nabla (\lambda) = H^0(\lambda)$. Moreover, as we shall
see (cf. Section 3.4 below) we have $\chi (\lambda) = \ch \Delta(\lambda) = \ch \nabla(\lambda)$.

\subsection{Chevalley groups}
Let $G_\Z$ be a split and connected reductive algebraic group scheme over $\Z$
corresponding to $G$. In other words $G_\Z$ is the associated Chevalley group.
Then $G$ is obtained from $G_\Z$ by extending scalars to $k$. More generally,
we write $G_A$ for the group scheme over an arbitrary commutative ring $A$
obtained via the base change $\Z \rightarrow A$. (The case $A = \Z_p$, the ring
of $p$-adic integers, will be needed in Chapter 5.) We use similar notation relative
to the subgroups $T$ and $B$. In particular, $T_\Z$ is a split maximal torus
in $G_\Z$ with $T_k = T$. We will identify $R$ with the root system associated
to $(G_A, T_A)$.

Note that for a $G_A$-module $V$ that is free of finite rank
as an $A$-module, $\ch(V)$ makes sense by considering ranks of weight spaces. If our 
field $k$ is an $A$-algebra then
we have for such a module $\ch(V) = \ch(V \otimes _A k)$.

For any commutative ring $A$ and any $B_A$-module $M$ we write $H^i_A(M)$ for
the $G_A$-module $R^i\Ind_{B_A}^{G_A} M$. See \cite{RAG}, I.5
for the general properties of these modules. In particular, we recall that if $A$ is
noetherian and $M$ is finitely generated over $A$, then $ H^i_A(M)$ is also
finitely generated over $A$, see \cite{RAG}, Proposition I.5.12 c).

Given any commutative ring $A$, for each $\lambda \in X^+$ we have the following
two $G_A$-modules: the Weyl module $\Delta_A(\lambda)$ and the dual Weyl
module $\nabla_A(\lambda)$. These modules are characteristic-free,
i.e., as $A$-modules both are free of rank equal to
$\dim \Delta(\lambda)$ and we have $G_A$-module isomorphisms
$\Delta_A(\lambda) \simeq \Delta_{\Z}(\lambda) \otimes A$
and $\nabla_A(\lambda) \simeq \nabla_{\Z}(\lambda) \otimes A$.
Just as for $G$, we have $\nabla_A (\lambda) = H_A^0(\lambda)$
and $\Delta_A (\lambda) \simeq H^N_A(w_0\cdot \lambda)$, see Chapter 3.

Any $G_\Z$-module $M$ which is finitely generated as a $\Z$-module has finite
torsion submodule $M_t = \{m\in M\mid nm = 0 \text { for some } n\in \N\}$. This
is a $G_\Z$-submodule and we set $M_f = M/M_t$. Then we refer to $M_t$ and
$M_f$ as the torsion part, respectively free part of $M$.

Any $M$ as above allows a surjection $P_0 \rightarrow M$ from a $G_\Z$-module
$P_0$ which is free of finite rank as a $\Z$-module. Hence $M$ also has a
free presentation $0 \rightarrow P_1 \rightarrow P_0 \rightarrow M \rightarrow 0$
with $P_0$ and $P_1$ free over $\Z$.

\subsection{Divisors}
Let $\mathcal  D(\Z)$ denote the divisor group for $\Z$, i.e., the free $\Z$-module with basis
consisting of all prime numbers $p$. If $n \in \Z\setminus {\{0\}}$ then we write
$\di (n) \in \mathcal  D(\Z)$ for the divisor associated to $n$. If $M$ is a finite
$\Z$-module of order $|M|$ we write $\di (M) = \di (|M|)$. Clearly,
$\di$ is additive with respect to short exact sequences of finite $\Z$-modules.

Suppose  $M$ is a  $T_\Z$-module. Then $M$ splits into a direct sum of weight submodules,
see \cite{RAG}, I.2.11.  When $M$ is finite, in analogy with the situation for $T$-modules 
in Section 2.2, this leads us to the following definition of $\di_T(M) \in \mathcal  D(\Z)[X]$
$$ \di_T(M) = \sum_{\mu \in X} \di(M_\mu) e^\mu.$$

Again it is clear that $\di_T$ is additive on exact sequences of finite
$T_\Z$-modules.

\subsection{Ext groups}

Consider finitely generated $G_{\Z}$-modules $M$ and $N$. By \cite {RAG},
II.B, the groups $\Ext_{G_\Z}^i(M,N)$ are finitely generated and vanish 
for large enough $i$. We will also need the following special cases of
some vanishing results from loc. cit.

\newpage

\begin{prop}
For $\lambda, \mu \in X^+$,
\begin{enumerate}
\item[a) ]  $\Ext^i_{G_\Z} (\Delta_\Z(\mu), \nabla_\Z(\lambda)) = 0$ 
unless ($\mu = \lambda$ and $i = 0$).
\newline $\Hom_{G_\Z} (\Delta_\Z(\lambda), \nabla_\Z(\lambda)) = \Z$.
\item[b) ] $\Ext^i_{G_\Z} (\Delta_\Z(\mu), \Delta_\Z(\lambda)) = 0$
  unless $\mu < \lambda$ or ($\mu = \lambda$ and $i = 0$).
  $\Hom_{G_\Z} (\Delta_\Z(\lambda),\allowbreak \Delta_\Z(\lambda)) = \Z$.
\end{enumerate}
\end{prop}

The universal coefficient theorem \cite{RAG}, Proposition I.4.18a 
gives analogous results over $G_A$ for other commutative rings $A$.
In particular the proposition stays valid after replacing each
$\Z$ by $\Z_p$. 

\subsection{Tilting modules}
A tilting module for $G_A$ is an $A$-finite $G_A$-module $Q$ which has both
a Weyl filtration (i.e., a filtration with successive quotients isomorphic to
Weyl modules) and a dual Weyl filtration (with successive quotients
isomorphic to dual Weyl modules). For a tilting module $Q$ (or more generally
any module with a Weyl filtration) we write $(Q:\Delta(\lambda))$ for the
number of times $\Delta_A(\lambda)$ occurs in a Weyl filtration of $Q$. This
integer is also uniquely defined by the character equation
$$ \ch Q = \sum_{\lambda \in X^+} (Q:\Delta(\lambda)) \chi (\lambda).$$

Let $A=\Z_p$ in this paragraph. We have the following standard facts, e.g., from
\cite{RAG}, II.E. For each $\lambda \in X^+$ there is a unique indecomposable
tilting module $T(\lambda)$ for~$G$ (respectively, $T_A(\lambda)$ for $G_A)$ with
highest weight $\lambda$. Every tilting module for~$G$ (respectively, of $G_A$)
is uniquely expressible as a direct sum of the various $T(\lambda)$ (respectively,
$T_A(\lambda)$). We have $T_A(\lambda) \otimes_A k \simeq T(\lambda)$. In particular,
every tilting module $\bar Q$ for $G$ lifts uniquely to a tilting module $Q$
for $G_A$ (i.e., $Q \otimes_A k \simeq {\bar Q}$).

\section{The Borel-Weil-Bott theorem and its consequences over $\Z$}

\subsection{The Borel-Weil-Bott theorem.}
The Borel-Weil-Bott theorem holds over any field of characteristic $0$. Here we
only state it over $\Q$. The general case then follows by an easy base change argument,
compare 3.2 below.
\begin{thm}[\cite{Bott}, \cite{Dem}]
  Let $\lambda \in X$ and choose $w \in W$ such that $w(\lambda +
  \rho) \in X^+$.  Then we have isomorphisms of $G_\Q$-modules
$$
H_\Q^i(\lambda) \simeq
\begin{cases}
H^0_\Q(w \cdot \lambda) &\text {if } i = l(w),\\
0 &\text { otherwise.}
\end{cases}
$$
\end{thm}
\begin{remark} Note that if $\lambda \in X$ is singular, i.e., if
  there exists $\alpha \in R$ with $\langle \lambda + \rho, \allowbreak
  \alpha^\vee \rangle =\nobreak 0$, then $H^i_\Q(\lambda) =0$ for all
  $i$.  Hence the possible non-uniqueness of $w$ in this statement
  does not cause ambiguity.
\end{remark}
\subsection{Universal coefficients theorem}
Let $A$ be an arbitrary commutative ring. Then for any $G_\Z$-module $M$ which is free of
finite rank over $\Z$ and for any $i \geq 0$ we have the following short exact
sequence of $A$-modules, cf. e.g., \cite{RAG} , I.4.18.
$$0 \rightarrow H^i_\Z(M) \otimes_\Z A \rightarrow H_A^i(M \otimes_\Z A)
\rightarrow \Tor_1^\Z (H^{i+1}_\Z(M), A) \rightarrow 0.$$

\subsection{The Borel-Weil-Bott theorem over $\Z$.}
When we combine 3.1 and 3.2 we find
\begin{cor}
Let $\lambda \in X$ and suppose $w \in W$ satisfies $w \cdot \lambda \in X^+$. Then
\begin{enumerate}
\item[a) ]
$H^i_\Z(\lambda)$ is a finite $\Z$-module for all $i \not = l(w)$.
\item[b) ]  $H^{l(w)}_\Z(\lambda)_f \otimes _\Z \Q \simeq \Delta_\Q(w\cdot \lambda)$.
\end{enumerate}
\end{cor}

\begin{rem}
If no $w \in W$ exists with $w \cdot \lambda \in X^+$ (i.e., if $\lambda$ is
singular) then $H^i_\Z(\lambda)$ is a finite $\Z$-module for all $i$.
\end{rem}

\subsection{Kempf's theorem}

Recall that Kempf's vanishing theorem \cite{Ke} says that if $\lambda$ is dominant then
all the higher cohomology modules $H^i(\lambda), i> 0$ vanish. This being true
for all fields we get (e.g., via the universal coefficient theorem above)

\begin{thm}
Let $\lambda \in X^+$. Then $H^i_\Z(\lambda) = 0$ for all $i>0$.
\end{thm}

This means in particular via the universal coefficient theorem above that for
dominant $\lambda$ we have that
$H^0_\Z(\lambda) \otimes_\Z k \simeq H^0(\lambda)$. Hence $\nabla_\Z(\lambda) \simeq
H^0_\Z(\lambda)$.

Via Serre duality (which is valid over all fields but not over $\Z$), Kempf's
theorem gives also
\begin{align}
 H^i(\lambda) = 0 \text { for all } i < N \text { and all } \lambda
\text { with } -\lambda -2 \rho \in X^+.
\end{align}
Hence for each $\lambda \in X^+$ we conclude that $H_\Z^N(w_0\cdot \lambda)$
has no torsion and the dimension of $H^N(w_0\cdot \lambda)$ is independent of
$k$. In fact, $H^N(w_0\cdot \lambda) \simeq H^N_\Z(w_0\cdot \lambda) \otimes _\Z k$
and hence  $H^N_\Z(w_0\cdot \lambda) \simeq \Delta_\Z(\lambda)$.

\begin{cor}
Let $V$ be any $G_\Z$-module (finitely generated over $\Z$ as always). Then for
all $\lambda \in X^+$ we have $H^i_\Z(V\otimes_\Z \lambda) = 0$  for all $i>0$ and
$H^0_\Z(V \otimes_\Z \lambda) \simeq V\otimes_\Z H^0_\Z(\lambda)$.
\end{cor}

\pf If $V$ is free over $\Z$ then we have the tensor identity \cite{RAG} I.3.6
$H^i_\Z(V\otimes_\Z \lambda) \simeq  V \otimes_\Z H^i_\Z(\lambda)$. Hence in this case
the corollary results directly from Kempf's theorem. In general, we have from
2.4 a presentation $0\rightarrow P_1 \rightarrow P_0 \rightarrow V \rightarrow0$
with $P_1$ and $P_0$ free over $\Z$. The corollary then holds for $P_1$ and $P_0$.
It is then immediate to deduce it for $V$.

Another important consequence of Kempf's theorem is that since it clearly gives
$\ch H^0(\lambda) = \sum_{i\geq 0} \ch H^i(\lambda) = \chi (\lambda)$ for all dominant
weights $\lambda$, we get
\begin{align}
\ch \nabla (\lambda) = \chi (\lambda) = \ch \Delta (\lambda) \text { for all } \lambda \in X^+.
\end{align}
Here the last equality follows by combining Kempf's vanishing and Serre duality,
see (3.1) above.

\begin{remark} It is well known that $\chi (\lambda)$ is given by Weyl's character formula, see, e.g.,
\cite {Do} (2.2.6).
\end{remark}

\subsection{Rank $1$}
The (very easy) proof by Demazure \cite{Dem} of Bott's theorem relies on an
analysis of natural isomorphisms $H^{i+1}_\Q(s_\alpha \cdot \lambda)
\rightarrow   H^{i}_\Q(\lambda)$ when $\alpha$ is a simple root with
$\langle \lambda, \alpha^\vee \rangle \geq 0$. We shall need the underlying
homomorphisms
over $\Z$ and hence engage in the following considerations.

Let $\alpha $ be a simple root and denote the corresponding minimal parabolic
subgroup $P_\alpha$ in $G$ containing $B$. Then we denote for any $B$-module
$M$ by $H^i_\alpha(M)$ the module $H^i(P_\alpha/B, M)$. Note that $P_\alpha/B$
is the projective line so that these cohomology modules always vanish for $i>1$.
When working over $\Z$ we write $H^i_{\alpha, \Z}(M)$ for the analogously
defined modules for the $\Z$-version $P_{\alpha,\Z}$ of $P_\alpha$.

\begin{lemma}[cf. \cite{RAG} II.5.2 and 8.13] Let $\lambda \in X$.
\begin{enumerate}
\item [a) ] If $\langle \lambda, \alpha^\vee \rangle \geq 0$ then
  $H^i_{\alpha, \Z} (\lambda) = 0$ for all $i >0$ and $H^0_{\alpha,
    \Z} (\lambda)$ is a free $\Z$-module whose weights are $\lambda,
  \lambda - \alpha, \dots , s_\alpha (\lambda)$, all occurring with
  multiplicity $1$.
\item [b) ] If $\langle \lambda, \alpha^\vee \rangle < -1$ then
  $H^i_{\alpha, \Z} (\lambda) = 0$ for all $i \not = 1$ and
  $H^1_{\alpha, \Z} (\lambda)$ is a free $\Z$-module whose weights are
  $\lambda + \alpha, \lambda + 2 \alpha, \dots , s_\alpha \cdot
  \lambda$, all occurring with multiplicity $1$.
\item [c) ] If $\langle \lambda, \alpha^\vee \rangle = r \geq 0$ then
  $\Hom_{P_{\alpha, \Z}}(H^1_{\alpha, \Z} (s_\alpha \cdot \lambda),
  H^0_{\alpha, \Z} (\lambda)) \simeq \Z$. Moreover, $H^0_{\alpha, \Z}
  (\lambda)$, respectively $H^1_{\alpha, \Z} (s_\alpha \cdot \lambda)$
  has a standard $\Z$-basis $\{v_0, v_1, \dots , v_r\}$, respectively
  $\{v_0', v_1', \dots , v_r'\}$ with $v_j$, respectively $v'_j$,
  having weight $\lambda - j \alpha$, $j = 0, 1, \dots, r$. A
  generator $c_\alpha(\lambda)$ of $\Hom_{P_{\alpha, \Z}}(H^1_{\alpha,
    \Z} (s_\alpha \cdot \lambda), H^0_{\alpha, \Z}(\lambda)) $ is
  given by
$$c_\alpha(\lambda)(v'_j) = \binom{r}{j} v_j, \ j= 0, 1, \dots , r.$$
\end{enumerate}
\end{lemma}

\subsection{Passing from Rank 1 to the general case.}
Keep the notation from 3.5. By transitivity of induction we have with obvious notation
$H^0(M) \simeq H^0(G/P_\alpha, H^0_\alpha(M))$. The same is true over $\Z$.
Hence using general properties of $H^i_\Z$, cf. \cite{RAG}, II.8 we
obtain from Lemma 3.6 a) and b)
that if $\lambda \in X$ and $\alpha \in S$ satisfy
$\langle \lambda, \alpha^\vee \rangle \geq 0$ then
\begin{align}
H^i_\Z(\lambda) \simeq H^i_\Z(H^0_{\alpha, \Z}(\lambda))
\end{align}
and
\begin{align}
H^{i+1}_\Z(s_\alpha \cdot \lambda) \simeq H^i_\Z(H^1_{\alpha, \Z}(s_\alpha \cdot\lambda))
\end{align}

Denote by $Q_\alpha(\lambda)$ the cokernel of the generator $c_\alpha(\lambda)$
from Lemma 3.6 c). Then $Q_\alpha(\lambda)$ is a finite $P_{\alpha, \Z}$-module
with weights $\lambda-\alpha, \lambda-2 \alpha, \dots , s_\alpha (\lambda)+\alpha$.
Each weight space is cyclic and we have
\begin{align}
\di _T(Q_\alpha(\lambda)) = \sum_{j =1}^{r-1} \di \binom{r}{j} e^{\lambda - j \alpha}.
\end{align}

The short exact sequence of $P_{\alpha, \Z}$-modules
$$ 0 \rightarrow H^1_{\alpha, \Z}(s_\alpha \cdot\lambda) \rightarrow
H^0_{\alpha, \Z}(\lambda) \rightarrow Q_\alpha(\lambda) \rightarrow 0$$
gives via (3.3) and (3.4) rise to the long exact sequence of $G_\Z$-modules
$$ \cdots \rightarrow H^{i+1}_\Z (s_\alpha \cdot \lambda) \rightarrow
 H^{i}_\Z (\lambda) \rightarrow  H^{i}_\Z (Q_\alpha(\lambda)) \rightarrow \cdots .$$

\begin{remark}
The isomorphisms  over $\Q$ analogous to (3.3) and (3.4) give isomorphisms
$H^{i+1}_\Q(s_\alpha \cdot \lambda) \simeq H^i_\Q(\lambda)$ for all $i$. This is
the key to Demazure's proof \cite{Dem} of Theorem~3.1.
\end{remark}

\section{Euler type formulas}

\subsection{Euler coefficients for $G$-modules.}

Let $V$ and $V'$ be $G_\Z$-modules, both finitely generated over $\Z$. Then
$\Ext_{G_\Z}^i(V, V')$ is finite for all $i>0$. This follows from Section 2.6
and the universal coefficient theorem \cite{RAG}, Proposition I.4.18a, because
$\Ext_{G_\C}^i(A, B) = 0$ for all $i>0$ and for any two rational $G_\C$-modules
$A$ and $B$ ($G_\C$ being reductive). If the $G_\C$-modules $V \otimes_\Z \C$
and $V' \otimes_\Z \C$ do not have an isomorphic simple summand, then
$\Hom_{G_\Z}(V, V')$ is finite. This happens in particular when $V$ or $V'$
is finite. By Section 2.6 we have in any case $\Ext^i_{G_\Z}(V, V') = 0$ when
$i \gg 0$. So whenever $\lambda \in X^+$ and $V$ is a $G_\Z$-module such that
$(V \otimes_\Z \C: \Delta_\C(\lambda)) = 0$ (e.g., when $V$ is a finite
$G_\Z$-module), the following expression gives a well defined element in $\Di (\Z)$
$$ e_\lambda^G(V) = \sum_{i\geq0} (-1)^i \di (\Ext_{G_\Z}^i(\Delta_\Z(\lambda), V)).$$
Clearly, $e_\lambda^G$ is additive on exact sequences of such $G_\Z$-modules
(in particular finite $G_\Z$-modules).

\begin{remark}
We may extend the above definition of $e_\lambda^G(V)$ to all (finitely generated)
$V$ by using just the torsion part of $\Hom_{G_\Z}(\Delta_\Z(\lambda), V)$. Clearly
when extended in this way $e_\lambda^G$ will fail to be additive on arbitrary exact 
sequences in general.
The proofs of Theorem 4.1 and Proposition 4.3 below require careful examination
of this failure for particular exact sequences.
\end{remark}

\subsection{Euler coefficients for $B$-modules.}
Suppose $M$ is a finite $B_\Z$-module. Then for each $i$ the $G_\Z$-module
$H_\Z^i(M)$ is also finite (because $G_\Z/B_\Z$ is a projective scheme) and it is $0$
for $i > N = \dim G_\Z/B_\Z$. We define for $\lambda \in X^+$
$$ e_\lambda^B(M) = \sum_{j\geq0}(-1)^j e^G_\lambda(H^j_\Z(M)).$$
Again we see that $e_\lambda^B$ is additive on exact sequences of finite $B_\Z$-modules.

If the $B_\Z$-structure on $M$ extends to $G_\Z$ then Corollary 3.5 tells
us that  $H^j_\Z(M) = 0$ for $j>0$ and
$H^0_\Z(M) \simeq M$. Hence in this case we have for all $\lambda \in X^+$
\begin{align}
e_\lambda^G(M) = e_\lambda^B(M).
\end{align}

\subsection{Formulas for $B$-modules.}
\begin{thm}
Let $M$ be a finite $B_\Z$-module. Then
\begin{enumerate}
\item [a) ] $e_\lambda^B(M) = \sum_{w \in W}(-1)^{l(w)} \di (M_{w\cdot \lambda})$
for all $\lambda \in X^+$.
\item [b) ]  $\sum_{i\geq0} (-1)^i \di _T(H^i_\Z(M))
= \sum_{\lambda \in X^+} e_\lambda^B(M) \chi(\lambda)$.
\end{enumerate}
\end{thm}

\pf a) The additivity of $e_\lambda^B$ immediately allows us to reduce to the case
where $M$ is given by the following short exact sequence
\begin{align}
0 \rightarrow \Z_\mu\stackrel{n} \rightarrow  \Z_\mu \rightarrow M \rightarrow 0
\end{align}
with $\mu \in X$ and $n\in \N$. In this case the formula we want to verify is
$$ e_\lambda^B(M) = \begin{cases} (-1)^{l(w)} \di (n) &\text{if } \mu = w\cdot \lambda
\text { for some } w\in W,\\
0& \text {otherwise.}\end{cases}$$
(Note that if $\mu = w\cdot\lambda$ for some $w\in W$ then $\mu$ is non-singular
and $w$ is uniquely determined. This is so because $\lambda \in X^+$).

To prove this we consider the long exact cohomology sequence arising from (4.2)
$$ \cdots \rightarrow H^i_\Z(\mu)\stackrel{n} \rightarrow H^i_\Z(\mu) \rightarrow H^i_\Z(M) \rightarrow \cdots.
$$
If $\mu$ is singular (i.e., if there is a $\beta \in R$ with
$\langle \mu + \rho, \beta^\vee \rangle = 0$) then all modules in this sequence are finite.
In this case the additivity of $e_\lambda^B$ immediately gives $e_\lambda^B(M) = 0$ as
desired.

So suppose $\mu$ is non-singular. Then there exists a unique $i_0$ such that
$H^{i_0}_\Z(\mu)$ is infinite. Define then $C^t(\mu)$, $C(\mu)$,
respectively $C^f(\mu)$ such that the diagram
\begin{center}
\parbox[t]{10cm}{
\xymatrix{&0\ar@{->}[d]& 0\ar@{->}[d]\\
&H^{i_0}_\Z(\mu)_t\ar@{->}[d]\ar@{->}[r]^{n}&H^{i_0}_\Z(\mu_{})_t\ar@{->}[d]\ar@{->}[r]
&C^t(\mu)\ar@{->}[d]\ar@{->}[r] &0\\
&H^{i_0}_\Z(\mu)\ar@{->}[d]\ar@{->}[r]^{n}&H^{i_0}_\Z(\mu_{})\ar@{->}[d]\ar@{->}[r]&
C(\mu)\ar@{->}[d]\ar@{->}[r] & 0\\
0\ar@{->}[r]&H_\Z^{i_0}(\mu)_f\ar@{->}[d]\ar@{->}[r]^{n}&H_\Z^{i_0}(\mu_{})_f\ar@{->}[d]\ar@{->}[r]&
C^f(\mu)\ar@{->}[r]& 0\\
&0 & 0}}
\end{center}
has exact rows. By definition the two first columns are exact and hence it follows that
so is the last column, i.e.,
$$0\rightarrow C^t(\mu) \rightarrow C(\mu) \rightarrow C^f(\mu) \rightarrow 0$$
is exact. Now the long exact sequence arising from (4.2) also gives exact sequences (recall
that $H^j_\Z(\mu)$ is a torsion module for $j \not = i_0$)
$$\cdots \rightarrow H^{i_0-1}_\Z(\mu) \rightarrow H^{i_0-1}_\Z(\mu)
\rightarrow H^{i_0-1}_\Z(M) \rightarrow H^{i_0}_\Z(\mu)_t \rightarrow
H^{i_0}_\Z(\mu)_t \rightarrow C^t(\mu) \rightarrow 0$$
and
$$0\rightarrow C(\mu) \rightarrow H^{i_0}_\Z(M) \rightarrow
H^{i_0+1}_\Z(\mu)\rightarrow H^{i_0+1}_\Z(\mu) \rightarrow H^{i_0+1}_\Z(M)
\rightarrow \cdots .$$
This implies
$$ e_\lambda^G(C^t(\mu)) = (-1)^{i_0-1} \sum_{j=0}^{i_0-1} (-1)^j e_\lambda^G(H^j_\Z(M))$$
and
$$e^G_\lambda(C(\mu)) =  (-1)^{i_0} \sum_{j\geq i_0} (-1)^j e_\lambda^G(H^j_\Z(M)).$$
We conclude that
$e_\lambda^B(M) = \sum_{j\geq 0} (-1)^j e_\lambda^G(H^j_\Z(M))=
(-1)^{i_0} (e^G_\lambda(C(\mu))-e_\lambda^G(C^t(\mu))) 
\allowbreak = (-1)^{i_0}e_\lambda^G(C^f(\mu))$.

So let $w \in W$ be determined by $w(\mu+\rho) \in X^+$. Then $i_0 = l(w)$ and
the weights of $H^{i_0}_\Z(\mu)_f$ coincide with those of $\Delta(w\cdot \mu)$.
In particular, if $\lambda = w\cdot \mu$ then $H^{i_0}_\Z(\mu)_f$ has unique
highest weight $\lambda$. Therefore
$\Hom_{G_\Z}(\Delta_\Z(\lambda),H^{i_0}_\Z(\mu)_f)
\simeq \Z$ and   $\Ext^i_{G_\Z}(\Delta_\Z(\lambda),H^{i_0}_\Z(\mu)_f) = 0$ for $i>0$
(since $\lambda \geq \nu$ for all weights $\nu$ of $H^{i_0}_\Z(\mu)_f$).
Hence in this case $e_\lambda^G(C^f(\mu)) = \di (n)$. If on the other hand
$\lambda \not = w\cdot \mu$ then $\Hom_{G_\Z}(\Delta_\Z(\lambda),H^{i_0}_\Z(\mu)_f)
= 0$ and the long exact $\Ext$-sequence arising from
$$0\rightarrow H^{i_0}_\Z(\mu)_f \stackrel{n}\rightarrow H^{i_0}_\Z(\mu)_f \rightarrow C^f(\mu) \rightarrow 0$$
consists entirely of
finite $\Z$-modules. It follows that in this case $e_\lambda^G(C^f(\mu))= 0$ as
desired.

b) Both sides of the equation in b) are additive in $M$ and hence just as
above we may restrict to the case where $M$ is defined by (4.2). Using a) for
the right hand side and recalling that $\chi (\mu) = (-1)^{l(w)}\chi(w\cdot \mu)$ for
all $w \in W$ we see that the desired equality in this case is
$$ \sum_{i\geq0}(-1)^i\di _T(H^i_\Z(M)) = \di (n) \chi(\mu).$$
Note that here $\chi (\mu) = 0$ if $\mu$ is singular.
Arguing as in the proof of a) above we obtain in fact
$$  \sum_{i\geq0}(-1)^i\di _T(H^i_\Z(M)) = (-1)^{i_0} \di _T(C^f(\mu)) =
\di (n) \chi(\mu).$$
Here the last equality results from the definition of $C^f(\mu)$ via
the fact that $\ch (H_\Z^{i_0}(\mu)_f) \allowbreak=
(-1)^{i_0} \chi (\mu)$.

\subsection{Formulas for $G$-modules.}
Combining Theorem 4.1 and the identity (4.1)
we obtain

\begin{cor} Let $V$ be a finite $G_\Z$-module
\begin{enumerate}
\item [a) ]
$e_\lambda^G(V) = \sum_{w \in W} (-1)^{l(w)}\di (V_{w\cdot \lambda})$ for all
$\lambda \in X^+$.

\item [b) ]
$\di _T(V) = \sum_{\lambda \in X^+} e_\lambda^G(V) \chi(\lambda)$.
\end{enumerate}
\end{cor}

\begin{remark}
The second identity in this corollary was obtained by the second author in \cite{Kul1}.
The argument there is different.
\end{remark}

\subsection{Natural homomorphisms.}
Fix now $\mu \in X^+$ and a reduced expression $w_0 = s_1s_2\cdots s_N$ for $w_0$ with
$s_i$ denoting the reflection corresponding to the simple root $\alpha_i$. Then we set
$$ \mu_0 = \mu, \ \mu_1=s_1\cdot \mu_0, \dots ,  \mu_i=s_i\cdot \mu_{i-1}, \dots ,
\mu_N=s_N\cdot \mu_{N-1} = w_0\cdot \mu.$$
Since $\mu$ is the unique highest weight of $\nabla_\Z(\mu)$ we have up to sign a unique
generator $c_\mu$ for $\Hom_{G_\Z}(\Delta_\Z(\mu), \nabla_\Z(\mu)) \simeq \Z$. We set
$Q(\mu) = \Coker(c_\mu)$ so that we have a short exact sequence
$$0\rightarrow \Delta_\Z(\mu) \rightarrow \nabla_\Z(\mu) \rightarrow Q(\mu) \rightarrow 0.$$
Now we claim that $c_\mu$ factors through $H_\Z^i(\mu_i)$ for all $i$. In fact, note that
$\langle \mu_{i-1}+\rho, \alpha_i^\vee\rangle =
\langle\mu+\rho, s_1s_2\cdots s_{i-1}(\alpha_i)^\vee\rangle > 0$ because
$s_1s_2\cdots s_{i-1}(\alpha_i) \in R^+$. Using the notation from Lemma 3.6 we therefore have
a short exact sequence
\begin{align}
0\rightarrow H_{\alpha_i, \Z}^1(\mu_i) 
\xrightarrow{ c_i }
H_{\alpha_i, \Z}^0(\mu_{i-1}) \rightarrow Q_{\alpha_i}(\mu_{i-1}) \rightarrow 0.
\end{align}
where $c_i = c_{\alpha_i}(\mu_{i-1})$ and $Q_{\alpha_i}(\mu_{i-1}) = \Coker (c_i)$. When we apply $H_\Z^{i-1}$
to (4.3) we get (see (3.3) and (3.4))
$$ 
\rightarrow H^i_\Z(\mu_i) 
\xrightarrow{ \tilde H^{i-1}_\Z(c_i)  }
H^{i-1}_\Z(\mu_{i-1})
\rightarrow H^{i-1}_\Z(Q_{\alpha_i}(\mu_{i-1})) \rightarrow 
$$
as part of a long exact sequence. Tracing a highest weight vector we see that (up to sign)
$c_\mu$ may be identified with the composite
$$ 
\Delta_\Z(\mu) \simeq H^N_\Z(\mu_N) 
\xrightarrow{\tilde c_N}
\cdots
\xrightarrow{\tilde c_{i+1}}
H^{i}_\Z(\mu_{i})_f
\xrightarrow{\tilde c_{i}}
H^{i-1}_\Z(\mu_{i-1})_f
\xrightarrow{\tilde c_{i-1}}
\cdots 
\xrightarrow{\tilde c_{i-1}}
H^0_\Z(\mu_0) \simeq \nabla_\Z(\mu).
$$
Note that we have passed to the free quotient of $H^i_\Z(\mu_i)$ and denoted the homomorphism
here induced by $H^{i-1}_\Z(c_i)$ by $\tilde c_i$. For $i=N$ and $i=0$
the cohomology modules are free, see Section 3.4 and so in these cases we have omitted the
subscript $f$. If $Q^f_i(\mu)$ denotes the cokernel of $\tilde c_i$ then we have a short exact sequence
$$ 0 \rightarrow  H^{i}_\Z(\mu_{i})_f \stackrel{\tilde c_{i}}\rightarrow H^{i-1}_\Z(\mu_{i-1})_f
\rightarrow Q^f_i(\mu) \rightarrow 0.$$

\subsection{Formulas for Euler coefficients.}
Keep the notation from 4.5. Then we have

\begin{prop}
$$  e^G_\lambda(Q(\mu)) = \sum_{i=1}^N (-1)^{i-1}e_\lambda^B(Q_{\alpha_i}(\mu_{i-1})) \text { for all }
\lambda \in X^+.$$
\end{prop}

\pf The factorization $c_\mu = \tilde c_1 \circ \tilde c_2 \circ \cdots \circ \tilde c_N$
from 4.5 gives immediately
\begin{align}
e^G_\lambda(Q(\mu)) = \sum_{i=1}^N e_\lambda^G(Q^f_i(\mu)) \text { for all }
\lambda \in X^+.
\end{align}

If we now set $Q_i(\mu) = \Coker(H^{i-1}_\Z(c_i))$ and let $Q_i^t(\mu)$ denote the
cokernel of the induced homomorphism
$H^i_\Z(\mu_i)_t \rightarrow H^{i-1}_\Z(\mu_{i-1})_t$ then we get the following commutative
diagram.
\begin{center}
\parbox[t]{10cm}{
\xymatrix{&0\ar@{->}[d]& 0\ar@{->}[d]& 0\ar@{->}[d]\\
&H^i(\mu_i)_t\ar@{->}[d]\ar@{->}[r]&H^{i-1}(\mu_{i-1})_t\ar@{->}[d]\ar@{->}[r]&
Q^t_i(\mu)\ar@{->}[d]\ar@{->}[r] &0\\
&H_\Z^i(\mu_i)\ar@{->}[d]\ar@{->}[r]&H_\Z^{i-1}(\mu_{i-1})\ar@{->}[d]\ar@{->}[r]&
Q_i(\mu)\ar@{->}[d]\ar@{->}[r] & 0\\
0\ar@{->}[r]&H^i(\mu_i)_f\ar@{->}[d]\ar@{->}[r]&H^{i-1}(\mu_{i-1})_f\ar@{->}[d]\ar@{->}[r]&Q^f_i(\mu)\ar@{->}[d]
\ar@{->}[r]& 0\\
&0 & 0&0}}
\end{center}
Here the rows and two first columns are exact. Hence we deduce that the last column
is also exact and we get
$$e^G_\lambda(Q_i^f(\mu)) = e^G_\lambda(Q_i(\mu)) -
e^G_\lambda(Q_i^t(\mu)).$$
Now exactly as in the proof of Theorem 4.1 the long exact sequences involved in the above diagram give
(note also that all the terms in the following expressions have to do with
finite $\Z$-modules)
\begin{align*}
e^G_\lambda(Q_i^t(\mu)) = {} & \sum_{j\leq i}
(-1)^{j-i-1}(e_\lambda^G(H^j_\Z(\mu_i)_t) -
e_\lambda^G(H^{j-1}_\Z(\mu_{i-1})_t))
\\ 
& +\sum_{j<i-1} (-1)^{j-i}e_\lambda^G(H_\Z^j(Q_{\alpha_i}(\mu_{i-1})))
\end{align*}
and
\begin{align*}
e^G_\lambda(Q_i(\mu)) = {} & \sum_{j>i}
(-1)^{j-i}(e_\lambda^G(H^j_\Z(\mu_i))
-e_\lambda^G(H^{j-1}_\Z(\mu_{i-1}))) 
\\
&- \sum_{j\geq i-1}
(-1)^{j-i}e_\lambda^G(H_\Z^j(Q_{\alpha_i}(\mu_{i-1}))).
\end{align*}
When we combine these two equations we obtain
\begin{align}
e^G_\lambda(Q_i^f(\mu)) = (-1)^i (e_i^t - e_{i-1}^t -e_\lambda^B(Q_{\alpha_i}(\mu_{i-1}))),
\end{align}
where we have set $e_r^t = \sum_{j\geq 0} (-1)^j e_\lambda^G(H^j_\Z(\mu_r)_t),\  r=0, 1, 2,
\dots ,N$.
When we now sum over $i$ in (4.5) we obtain the desired equality since the
$e^t_i$ cancel each other
leaving only $e^t_0$ and $e^t_N$. Both these are $0$ by Kempf's theorem, see Section 3.4.

\subsection{Two lemmas.}
We still use the notation from 4.5. Now we shall combine Proposition 4.3 with Theorem 4.1.
Recall from (3.5) that the weights of $Q_{\alpha_i}(\mu_{i-1})$ are
$\mu_{i-1} - \alpha_i, \mu_{i-1} - 2\alpha_i, \dots , \mu_{i-1} - (r_i-1)\alpha_i$ where
$r_i = \langle\mu_{i-1}, \alpha_i^\vee\rangle$. All weight spaces are
cyclic and the order of $Q_{\alpha_i}(\mu_i)_{\mu_{i-1} - m\alpha_i}$ is
$\binom{r_i}{m}$, $m = 1, 2, \dots , r_i-1$.

\begin{lemma}
Let $\lambda \in X^+$ and  $x, y \in W$. Suppose that both $x\cdot \lambda$ and $y\cdot \lambda$ are weights
of $Q_{\alpha_i}(\mu_{i-1})$. Then either $x=y$ or $x = s_iy$.
\end{lemma}

\pf Suppose $x\cdot \lambda = \mu_{i-1} - m\alpha_i$ and
$y\cdot \lambda = \mu_{i-1} - m'\alpha_i$ with $0<m, m'<r_i$. Then
$x\cdot \lambda = y\cdot \lambda + (m'-m)\alpha_i$. Hence
$$
(\lambda+\rho,\lambda + \rho) = (\lambda+\rho, \lambda + \rho) + (\alpha_i,\alpha_i)
(m'-m)(\langle y(\lambda+ \rho),  \alpha_i^\vee
\rangle +(m'-m))
$$ 
and we conclude 
that either $m'=m$  or $m-m' =
\langle y(\lambda+\rho), \alpha_i^\vee \rangle$. In the first case $y\cdot \lambda
= x\cdot \lambda$ and therefore $y = x$. In the second case we get
$s_iy\cdot \lambda =
y\cdot \lambda - \langle y(\lambda+\rho), \allowbreak \alpha_i^\vee \rangle \alpha_i =
y\cdot \lambda -(m-m')\alpha_i = x\cdot \lambda$, i.e., $y= s_ix$.

\begin{lemma} Let $\lambda \in X^+$. Suppose there exist $x \in W$ and $0< m < r_i$ with
$x\cdot \lambda = \mu_{i-1} -m\alpha_i$. Then
 \begin{eqnarray*}
\sum_{w\in W} (-1)^{l(w)}\di(Q_{\alpha_i}(\mu_{i-1})_{w\cdot\lambda}) =
      (-1)^{l(x)}(\di(r_i+1-m) - \di(m)).
\end{eqnarray*}
\end{lemma}

\pf This follows from the Lemma 4.4 together with the observation that for all $r \geq m \geq 0$
we have
$$ \di (\binom{r}{m}) - \di (\binom{r}{r+1-m}) = \di (r+1-m) - \di (m).$$
Note in particular that the lemma holds also when $m =1$ (in which case
$s_ix\cdot\lambda$ is not a weight of $Q_{\alpha_i}(\mu_i)$).

\subsection{An Euler type formula.}
Let $\lambda, \mu \in X^+$. For each $\beta \in R^+$ we set
$$ V_\beta(\lambda, \mu) = \{(x,m)\mid x\in W, 0<m< \langle\mu+\rho,\beta^\vee\rangle \text { with }
x\cdot \lambda = \mu - m\beta \}.$$

With this notation we have
\begin{thm}
The cokernel $Q(\mu)$ of the canonical homomorphism $\Delta_\Z(\mu) \rightarrow \nabla_\Z(\mu)$
satisfies
$$e_\lambda^G(Q(\mu)) = -\sum_{\beta \in R^+} \sum_{(x,m) \in V_\beta(\lambda, \mu)} (-1)^{l(x)}\di (m).$$
\end{thm}

\pf
When we combine Theorem 4.1 a) and Proposition 4.3 we get
$$ e^G_\lambda(Q(\mu)) = \sum_{i=1}^N(-1)^{i-1}
\sum_{w\in W}(-1)^{l(w)} \di (Q_{\alpha_i}(\mu_{i-1})_{w\cdot \lambda}).$$
Note that if we set $\beta_i = s_1s_2\cdots s_{i-1}(\alpha_i)$ then
$\{\beta_1, \beta_2, \dots ,\beta_n\} = R^+$. Moreover, the equality
$x\cdot \lambda =
\mu_{i-1} - m\alpha_i$
is equivalent to $s_1s_2\cdots s_{i-1}x\cdot \lambda = \mu - m\beta_i$. Also $r_i =
\langle \mu_{i-1},\alpha_i^\vee \rangle = \langle(\mu+\rho,\beta_i^\vee \rangle -1$. Hence the
theorem follows by Lemma 4.4 and Lemma 4.5.

\begin{rem}
The arguments in Lemma 4.4 show that the set $V_\beta(\lambda, \mu)$ is either
empty or contains exactly two elements (of the form $(x,m)$ and
$(s_\beta x, \langle \mu + \rho, \beta^\vee \rangle - \nobreak m)$).
\end{rem}

\subsection{Variations.} We present some variations of Theorem 4.6 for later use.

\begin{cor}
$$e_\lambda^G(\Delta_\Z(\mu)) = \sum_{\beta \in R^+} \sum_{(x,m) \in V_\beta(\lambda, \mu)} (-1)^{l(x)}\di (m).$$
\end{cor}

\pf
Use Proposition 2.1 with the sequence
$0 \rightarrow \Delta_\Z(\mu) \stackrel{c_\mu} \rightarrow \nabla_\Z(\mu) \rightarrow Q(\mu) \rightarrow 0$.
(Note that the corollary--as understood by the Remark in Section 4.1--and its proof
are valid even for $\lambda = \mu$. We have $e_\lambda^G(\Delta_\Z(\lambda)) = 0$,
see Proposition 7.1 below.)

Let $\lambda, \mu \in X^+$. For each $\gamma \in R^+$ we set
$$ U_\gamma(\lambda, \mu) = \{(w,n)\mid w\in W, n < 0 \text { or }
n > \langle\lambda+\rho,\gamma^\vee\rangle, w\cdot \mu = \lambda - n\gamma \}.$$
With this notation, we can deduce from Theorem 4.6 an alternate expression for
$e_\lambda^G(\Delta_\Z(\mu))$.

\begin{prop}
$$e_\lambda^G(\Delta_\Z(\mu)) = \sum_{\gamma \in R^+} \sum_{(w,n) \in U_\gamma(\lambda, \mu)} (-1)^{l(w)}\di (n).$$
\end{prop}

\pf Let
$V(\lambda, \mu) = \bigcup_{\beta \in R^+} \{(\beta,x,m)\mid (x,m) \in
V_\beta(\lambda, \mu) \}$ 
and
\[
U(\lambda, \mu) = \bigcup_{\gamma \in R^+} \{(\gamma,w,n)\mid (w,n)
\in U_\gamma(\lambda, \mu) \}.
\] 
By Corollary 4.8 it is enough to produce a bijection between $U(\lambda,\mu)$ and
$V(\lambda,\mu)$ for which $m = \pm n$ and $x = w^{-1}$. This is an easy check as
follows.

First let $(\gamma,w,n) \in U(\lambda,\mu)$. Since $\lambda - n \gamma = w \cdot \mu$,
we have $w^{-1} \cdot \lambda = \mu + n(w^{-1}\gamma)$.

{\it Case 1a.} If $w^{-1}\gamma \in R^+$ then let $\beta = w^{-1}\gamma$, $x = w^{-1}$
and $m = -n$. We have
\begin{align}
\langle \mu + \rho , w^{-1}\gamma^\vee \rangle =
\langle w^{-1}(\lambda + \rho) - n(w^{-1}\gamma), w^{-1}\gamma^\vee\rangle =
\langle \lambda + \rho, \gamma^\vee\rangle - 2n.
\end{align}
Since $\langle \mu+\rho , \beta^\vee \rangle > 0$ and $\langle \lambda+\rho, \gamma^\vee\rangle > 0$,
the possibility $n > \langle \lambda + \rho, \gamma^\vee\rangle$ in the definition of
$U_\gamma(\lambda, \mu)$ cannot be true. So $n < 0$ and hence $m = -n >0$. Also
$\langle \mu + \rho , \beta^\vee\rangle = \langle \lambda + \rho, \gamma^\vee\rangle - 2n > -2n = 2m$.
So $0 < m < \frac{1}{2} \langle \mu + \rho , \beta^\vee\rangle$; in particular
$(\beta,x,m) \in V(\lambda,\mu)$.

{\it Case 1b.} If $w^{-1}\gamma \in -R^+$ then let $\beta = -w^{-1}\gamma$, $x = w^{-1}$
and $m = n$. By (4.6) we have
$\langle \mu + \rho , \beta^\vee \rangle = 2n - \langle \lambda + \rho, \gamma^\vee \rangle.$ Since
$\langle \mu + \rho , \beta^\vee \rangle > 0$ and  $\langle \lambda + \rho, \gamma^\vee\rangle > 0$,
the possibility $n < 0$ in the definition of $U_\gamma(\lambda, \mu)$ cannot be true.
So $n > \langle \lambda + \rho, \gamma^\vee \rangle$. Thus $m = n > 0$. Also
$\langle \mu + \rho , \beta^\vee \rangle = 2n - \langle \lambda + \rho, \gamma^\vee \rangle > n = m$,
as desired. (Since $0< \langle \lambda +\rho,\gamma^\vee \rangle = 2n- \langle \mu+\rho,\beta^\vee \rangle$,
we actually have
$\frac{1}{2} \langle \mu +\rho , \beta^\vee \rangle < n = m < \langle \mu + \rho, \beta^\vee \rangle$.)

For the inverse map, let $(\beta, x, m) \in V(\lambda,\mu)$. Since $\mu - m\beta = x \cdot \lambda$,
we have $x^{-1} \cdot \mu = \lambda + m(x^{-1}\beta)$.

{\it Case 2a.} If $x^{-1}\beta \in R^+$ then let $\gamma = x^{-1}\beta$, $w = x^{-1}$ and
$n = -m$. We have $n < 0$ since $m > 0$, so $(\gamma,w,n) \in U(\lambda ,\mu)$. Clearly
this case is inverse to Case 1a. ($m$ must satisfy the bounds in the last sentence of Case
1a by the calculation there.)

{\it Case 2b.} If $x^{-1}\beta \in -R^+$ then let $\gamma = -x^{-1}\beta$, $w = x^{-1}$
and $n = m$. Now via a calculation similar to (4.6) we have
$\langle \lambda + \rho, \gamma^\vee\rangle = 2m - \langle \mu + \rho , \beta^\vee\rangle < 2m-m = n$.
so $(\gamma,w,n) \in U(\lambda ,\mu)$. Clearly this case is inverse to Case 1b (and again
the bounds obtained there on $m$ must hold in this case).

\begin{remark}
Note that the above bijection pairs $V_\beta(\lambda, \mu)$ and
$U_\gamma(\lambda, \mu)$ where $\gamma = x^{-1}\beta$ with $x \in W$
chosen such that $(x,m) \in V_\beta(\lambda, \mu)$ and  $x^{-1}\beta
\in R^+$. (This is always possible by replacing $x$ with $s_\beta x$ if
necessary, see Remark 4.7.)
\end{remark}

\section{Sum Formulas}

\subsection{Sum formula for Weyl modules}
Let $\Delta^i_\Z(\mu) = c_\mu^{-1} (p^i \nabla_\Z(\mu)) \subset \Delta_\Z(\mu)$.
Jantzen's filtration is a descending filtration of $\Delta(\mu)$
defined by $\Delta^i(\mu)=$ the $G$-submodule generated by the
image of $\Delta^i_\Z(\mu)$ under the canonical map
$\Delta_\Z(\mu) \rightarrow \Delta(\mu)$. We now have Jantzen's
sum formula cf. \cite{JCJ} and \cite{Asum}.

\begin{cor} Let $\nu_p$ denote the $p$-adic valuation. Then
$$
\sum_{i>0} \ch(\Delta^i(\mu)) = \sum_{\beta \in R^+}
\sum_{0 < m < \langle \mu + \rho, \beta^\vee \rangle}
\nu_p(m) \chi(\mu - m \beta).
$$
\end{cor}

\pf
It is well-known that the left hand side is the coefficient of
$[p]$ in $\di_T(Q(\mu))$, e.g., diagonalize $c_\mu$ and calculate
each expression. The result follows by Corollary 4.2b
and Theorem 4.6.

\subsection{A filtration associated to tilting modules}

Henceforth in this chapter  we let $A = \Z_p$, the ring of
$p$-adic integers. Fix $Q$, a tilting $G_A$-module.
Also fix $\lambda \in X^+$. Following \cite{Atilt}
we define two descending filtrations as follows.
First, let $F_{\lambda}(Q) = \Hom_{G_A}(\Delta_A(\lambda), Q)$.
Define
$$F_{\lambda}(Q)^j = \{\varphi \in F_{\lambda}(Q) \mid
\psi \circ \varphi \in p^j A c_\lambda \text{ for all }
\psi \in \Hom_{G_A}(Q,\nabla_A(\lambda))\},$$
where $c_\lambda$ now denotes a generator of
$\Hom_{G_A}(\Delta_A(\lambda), \nabla_A(\lambda))$.
Next, let ${\bar Q} = Q \otimes_A k$. Recall that $Q$
is determined uniquely by ${\bar Q}$. Let
${\bar F}_{\lambda}({\bar Q})=\Hom_G(\Delta(\lambda),{\bar Q})
=\Hom_{G_A}(\Delta_A(\lambda), Q) \otimes_A k$ (by, e.g.,
universal coefficients and Proposition 2.1 a)).
Define ${\bar F}_{\lambda}({\bar Q})^j =$ the $k$-vector
space spanned by the image of $F_{\lambda}(Q)^j$ in
$F_{\lambda}(Q)\otimes_A k$. In the remaining sections
we will prove a sum formula for the latter filtration.

\subsection{Homological considerations}
Continue with the notation in 5.2. Following~\cite{Atiltsum}, Chapter
1, we relate the desired sum formula to certain Ext groups via an
equivalent description of $F_{\lambda}(Q)^j$. For this, fix an
enumeration of dominant weights such that $\lambda_i < \lambda_j$
implies $i < j$. Let $(Q:\Delta_A(\lambda_j)) = n_j$.  We will freely
use Proposition 2.1 in the following analysis without further mention.
A first application gives that $Q$ has a finite filtration $Q = Q_0
\supset Q_1 \supset Q_2 \cdots$ with $Q_{i-1}/Q_i =
\Delta_A(\lambda_i)^{n_i}$ for some $n_i \geq 0$.  Now fix $i$ such
that the chosen $\lambda=\lambda_i$. Consider the two short exact
sequences
\begin{align}
0 \rightarrow Q_{i-1} \rightarrow Q \rightarrow Q/Q_{i-1} \rightarrow 0
\quad \text{ and } \quad
0 \rightarrow Q_i \rightarrow Q_{i-1} \stackrel{\pi} \rightarrow
\Delta_A(\lambda)^{n_i} \rightarrow 0.
\end{align}
Apply $\Hom_{G_A}(\Delta_A(\lambda), -)$ to these.
In the first long exact sequence, for $t > 0$,
$\Ext^t_{G_A}(\Delta_A(\lambda), Q/Q_{i-1}) = 0$ (since
$(Q/Q_{i-1}: \Delta_A(\lambda_j)) = 0$ for any
$\lambda_j > \lambda$) and $\Ext^t_{G_A}(\Delta_A(\lambda), Q) = 0$
(since $Q$ has a dual Weyl filtration). Also
$\Hom_{G_A}(\Delta_A(\lambda), \allowbreak Q/Q_{i-1}) = 0$, since
$(Q/Q_{i-1}: \Delta_A(\lambda)) = 0$. Hence
$\Ext^t_{G_A}(\Delta_A(\lambda), Q_{i-1})=0$ for $t>0$
and the entire sequence reduces to the isomorphism
$\Hom_{G_A}(\Delta_A(\lambda),Q_{i-1})\simeq F_{\lambda}(Q)$.

Next we use this information in the second long exact
sequence. Since $\Ext^t_{G_A}(\Delta_A(\lambda), \allowbreak \Delta_A(\lambda)) = 0$
for $t>0$, we get $\Ext^t_{G_A}(\Delta_A(\lambda), Q_i) = 0$
for $t>1$. Also $\Hom_{G_A}(\Delta_A(\lambda), \allowbreak Q_i) = 0$
since $(Q_i: \Delta_A(\lambda)) = 0$. So the entire
sequence reduces to
\begin{align}
0 \rightarrow F_{\lambda}(Q) \stackrel{\Phi} \rightarrow
\End_{G_A}(\Delta_A(\lambda))^{\oplus{n_i}}
\rightarrow \Ext^1_{G_A}(\Delta_A(\lambda), Q_i) \rightarrow 0.
\end{align}
Still following \cite{Atiltsum}, we take a closer look at
certain maps between several Hom-groups. First, note that
$\Phi(\varphi) = \pi \circ \varphi$. (This makes
sense  since, by the previous paragraph, any map
$\varphi \in F_{\lambda}(Q)$ factors through $Q_{i-1}$.)
Note that $\Phi$ is an injection between free $A$-modules,
each of rank $n_i$. Next, apply $\Hom(-,\nabla_A(\lambda))$
to the two short exact sequences (5.1). Each of the
resulting long exact sequences reduces to just Hom-terms.
Since $(Q_i:\Delta_A(\lambda)) = 0 = (Q/Q_{i-1}:\Delta_A(\lambda))$,
we get isomorphisms
$\Hom_{G_A}(Q,\nabla_A(\lambda)) \simeq
\Hom_{G_A}(Q_{i-1},\nabla_A(\lambda)) \simeq
\Hom_{G_A}(\Delta_A(\lambda),\nabla_A(\lambda))^{\oplus n_i}.$
This sequence of bijections pairs $\psi \in \Hom_{G_A}(Q,\nabla_A(\lambda))$
first with its restriction $\psi |_{Q_{i-1}}$ and then to
${\bar \psi} \in \Hom_{G_A}(\Delta_A(\lambda)^{\oplus n_i},\nabla_A(\lambda))$
such that ${\bar \psi} \circ \pi = \psi |_{Q_{i-1}}$.
So $\psi \circ \varphi = {\bar \psi} \circ \pi \circ \varphi = {\bar \psi} \circ \Phi(\varphi)$.
This easily gives (see \cite{Atiltsum}, Proposition 1.6):
\begin{align}
F_{\lambda}(Q)^j = \{ \varphi \in F_{\lambda}(Q) \mid \Phi(\varphi)
\in p^j \End_{G_A}(\Delta_A(\lambda))^{\oplus {n_i}} \}.
\end{align}

\subsection{A sum formula involving tilting modules}
Keep the notation from 5.2 and 5.3. Additionally, for
arbitrary $\xi \in X$, we make the following notation.
If there exists $w \in W$ with $\mu = w \cdot \xi$ dominant,
define $[Q: \chi(\xi)] = (-1)^{\ell(w)} (Q: \Delta_A(\mu))$.
Otherwise let $[Q: \chi(\xi)] = 0$. This makes sense by
Theorem 3.1. We now prove the following sum formula,
which was discovered (and proved when $p \geq h$) in \cite{Atiltsum}.
\begin{thm}
$$
\sum_{j > 0} \dim  {\bar F}_{\lambda}({\bar Q})^j =
- \sum_{\alpha \in R^+}
\sum_{n < 0 \text{ \rm or } n > \langle \lambda + \rho, \alpha^\vee \rangle}
\nu_p (n) [{\bar Q} : \chi(\lambda - n \alpha)].
$$
\end{thm}
\pf
From (5.2) and (5.3) it is standard (e.g., by diagonalizing $\Phi$) to see that
$$\sum_{j > 0} \dim {\bar F}_{\lambda}({\bar Q})^j =
\nu_p(\Ext^1_{G_A}(\Delta_A(\lambda), Q_i)).$$
Since $\Ext^t_{G_A}(\Delta_A(\lambda), Q_i) = 0$ for $t \neq 1$,
we have
$$\nu_p(\Ext^1_{G_A}(\Delta_A(\lambda), Q_i))
= - \sum_t (-1)^t \nu_p(\Ext^t_{G_A}(\Delta_A(\lambda), Q_i)).$$
Recall that $(Q_i:\Delta(\lambda_j))$ is $n_j$ if $j>i$ and 0
otherwise. So
$$\sum_t (-1)^t \nu_p(\Ext^t_{G_A}(\Delta_A(\lambda), Q_i))
= \sum_{j>i} n_j \sum_t (-1)^t \nu_p(\Ext^t_{G_A}(\Delta_A(\lambda), \Delta_A(\lambda_j))).$$
(Note that all the Hom-terms in the previous equation are zero,
so additivity of Euler characteristic holds.) The last alternating sum
in the preceding equation may be replaced by the coefficient
of $[p]$ in $e_\lambda^G(\Delta_\Z(\lambda_j))$.
Then we may take the outer sum over
all $j$ as $e_\lambda^G(\Delta_\Z(\lambda_j)) = 0$ for $j \leq i$.
Altogether we have
$$\sum_{j > 0} \dim  {\bar F}_{\lambda}({\bar Q})^j = -
\text{ the coefficient of } [p] \text{ in }
\sum_j n_j e_\lambda^G(\Delta_\Z(\lambda_j)).$$
The result follows by Proposition 4.9.

\section{Quantum Groups}

\subsection{Passing to the quantum case.}
The sum formulas Corollary 5.1 and Theorem~5.2 have direct analogues
for quantum groups at roots of $1$, see \cite{APW}, \cite{T} and
\cite{Atiltsum}. We shall show in this section that our approach above
carries over to the quantum case. In particular, this allows us to get
rid of the condition in loc. cit. that the order of the root of unity
must be at least equal to the Coxeter number. In the Weyl module case
the reason for this restriction was that the quantized Kempf vanishing
theorem had only been proved in that case (see \cite{APW} and
\cite{AW}). This restriction was removed by Ryom-Hansen's general
proof \cite{RH}. In the tilting module case the reason for the
restriction was that the proof in \cite{Atiltsum} required a regular
weight as its starting point.

We carefully set up the quantized version of the approach in Sections 6.2--5. Once
this is done the arguments are completely parallel and we shall leave to the reader
the task of repeating the proofs leading to the quantized versions of the sum
formulas.

\subsection{The quantum parameter.}
Throughout this section $k$ will denote an arbitrary field. We set $p = \cha (k)
\geq 0$. For technical reasons we need $p \not = 2$ and also that $p
\neq 3$ if the root system in question contains type $G_2$. Then we fix a root of
unity $q \in k$ of order $l$ or $2l$ with $l \in \N$ odd.

We let $v$ denote an indeterminate and set $\mathcal A = \Z[v,v^{-1}], \; A=k[v,v^{-1}]$.
The natural homomorphism $\mathcal A \rightarrow A$ mapping $v \in
\mathcal A$ to $v \in A$ makes $A$ into an $\mathcal A$-algebra. We make $k$ into
an $\mathcal A$-algebra by specializing $v$ to $q$. Of course so far $q$ could be
any non-zero element in $k$ but as we shall see the only interesting case for our
present purposes is when $q$ is a root of unity.

\subsection{Roots and weights}
As in Section 2.1 we denote by $R$ a (finite) root system and we choose a set of
positive roots $R^+$. This takes place in some euclidian space $E = \R^n$ and
we let $\{\alpha_1, \alpha_2, \dots , \alpha_n\}$ be an enumeration of the set
of simple roots $S \subset R^+$. Moreover, we denote by $X \subset E$ the set of
integral weights, i.e.,
$$ X = \{\lambda \in E \mid \langle \lambda, \beta^\vee \rangle \in \Z,\  \beta \in R\}.$$
Then $X \simeq \Z^n$. As before we set $X^+$ equal to the set of dominant
weights in $X$.

The Weyl group $W$ of $R$ acts naturally on $E$ and $X$. Again we also have the dot-action
given by $w\cdot\lambda = w(\lambda + \rho) - \rho, \ w\in W, \lambda \in E$
with $\rho = \frac{1}{2} \sum_{\beta\in R^+} \beta$.

\subsection{Quantum groups over $k$.}
Let $U$ denote the quantum group over $\Q(v)$ associated with $R$. This is the
$\Q(v)$-algebra defined by some generators $E_i, F_i, K_i^{\pm}, \
i= 1, 2, \dots , n$ and certain relations, see e.g. \cite{L}. It has a triangular
decomposition $U = U^-U^0U^+$ with $U^-$, respectively $U^0$, $U^+$ denoting  the
subalgebra generated by all $F_i$'s, respectively $K_i^{\pm}$'s,  $E_i$'s.

Inside $U$ we have an $\mathcal A$-subalgebra $U_\mathcal A$, the Lusztig
$\mathcal A$-form of $U$. It is defined via the (gaussian) divided powers
$E_i^{(m)}$ and $F_i^{(m)}$, $m \in \N,\; i= 1, 2, \dots , n$, see \cite{L}. Then
for each $\mathcal A$-algebra $A'$ we set $U_{A'} = U_\mathcal A\otimes_\mathcal AA'$
and call this the quantum group over $A'$ associated with $R$. In particular, $U = U_{\Q(v)}$
when $\Q(v)$ is given the natural $\mathcal A$-structure obtained by sending
$v$ to $v$. In the case where $A' = k$ with $\mathcal A$-structure as above
we often write $U_q$ instead of $U_k$.

The above triangular decomposition of $U$ generalizes to $U_{A'} =
U_{A'}^-U_{A'}^0U_{A'}^+$
for appropriate $A'$-subalgebras $U_{A'}^-, U_{A'}^0$, and $U_{A'}^+$, see
\cite{L}. We set $B_{A'} = U^-_{A'}U^0_{A'}$. Again we write $B_q$ instead
of $B_k$.

\subsection{Integrable modules and induction functors.}
Let $\lambda \in X$. Then for any $A'$ as above $\lambda$ gives rise to a character
$\chi_\lambda : U_{A'}^0 \rightarrow A'$ which extends uniquely to a character of
$B_{A'}$ (taking all $F_i^{(m)}$'s to $0$). Then
if $M$ is a $U_{A'}^0$-module we set
$$ M_\lambda = \{m\in M \mid um = \chi_\lambda(u) m, \ u\in U_{A'}^0 \}$$
and call this the $\lambda$-weight space in $M$.

We denote by $\mathcal C_{A'}$, respectively $\mathcal C^-_{A'}$ the category
consisting of all integrable $U_{A'}$-, respectively $B_{A'}$-modules. A module
is integrable if it splits into a direct sum of its weight spaces (as $U^0_{A'}$-module)
and all high enough divided powers of all relevant generators vanish on any given
element in the module, see \cite{APW}.

We have a natural induction functor $\Ind_{B_{A'}}^{U_{A'}}: \mathcal C_{A'}^-
\rightarrow \mathcal C_{A'}$, see \cite {APW}. As in loc. cit. we shall denote the
right derived functors of this functor by  \smash{$H^j_{A'}$}, $j \geq 0$. These functors
share many of the properties of the $G/B$-sheaf cohomology functors from Section~2.3.
In particular, we have (cf. \cite{APW} Theorem 5.8)
\begin{enumerate}
\item
If $M \in \mathcal C^-_{A'}$ is finitely generated as an $A'$-module then each
$H^j_{A'}(M) \in \mathcal C_{A'}$ is also finitely generated over $A'$,
\item
$H^j_{A'} = 0$ for all $j > N$
\end{enumerate} (as before $N$ denotes the number of positive roots).

\subsection{Weyl and dual Weyl modules.}
Keep the notation from above and fix now $\mu \in X^+$. Then we set
$$ \Delta_{A'} (\mu) = H^N_{A'}(w_0\cdot \mu) \text { and } \nabla_{A'}(\mu) =
H^0_{A'}(\mu).$$
We call these the Weyl module and the dual Weyl module for $U_{A'}$ with highest weight
$\mu$. Because of the
quantized Kempf's vanishing theorem (which was proved for special $A'$'s in \cite {APW}
and in general by Ryom-Hansen in \cite{RH})we have in analogy with Section 3.4
\begin{align}
\nabla_{A'}(\mu) = \nabla_\mathcal A(\mu) \otimes_\mathcal A A'.
\end{align}
Since $H^{N+1}_{\mathcal A} = 0$ we also have
\begin{align}
\Delta_{A'}(\mu) = \Delta_\mathcal A (\mu) \otimes_\mathcal A A'.
\end{align}
Just as in Section 3.4 we get quantized Weyl character formulas
\begin{align}
\ch \Delta_q (\mu) = \ch \nabla_q (\mu) = \chi (\lambda).
\end{align}

\subsection{Simple modules.}
We shall now consider the case $A' = k$. Then each $\nabla_q(\mu)$ (again we
use index $q$ instead of $k$ here) contains a unique simple submodule which
we denote $L_q(\mu)$, cf. \cite{APW}. The family $\{L_q(\mu)\}_{\mu \in X^+}$
is then up to isomorphisms the set of simple modules in $\mathcal C_q$
(and in fact this family and their sign-twists constitute all finite dimensional
simple $U_q$-modules, see loc. cit.).

Serre duality gives that $L_q(\mu)$ is also the unique simple quotient of
$\Delta_q(\mu)$. In fact, we have up to scalars a unique homomorphism $c_q(\mu) :
\Delta_q(\mu) \rightarrow \nabla_q(\mu)$ and the image is $L_q(\mu)$. This
homomorphism may be obtained by specialization from a generator
$c_{\mathcal A} (\mu)$ of
$\Hom_{\mathcal C_\mathcal A}(\Delta_\mathcal A(\mu), \nabla_\mathcal A(\mu))$. We
shall now study the corresponding homomorphism $c_\mu = c_A(\mu) \in
\Hom_{\mathcal C_A}(\Delta_A(\mu), \nabla_A(\mu))$ (obtained from $c_\mathcal A(\mu)$
by the base change $\mathcal A \rightarrow A$) just as we studied the corresponding
homomorphism in Section~4.5.

\begin{remark}
If $q$ was not a root of unity then we would have $\nabla_q(\mu) = L_q(\mu) =
\Delta _q(\mu)$ for all $\mu \in X^+$, and $\mathcal C_q$ would be semisimple,
see \cite{APW}.
\end{remark}

\subsection{Rank 1.}
Fix $i \in \{1, 2, \dots , n\}$. Then we set $P_{A}(i)$ equal to the $A$-subalgebra
of $U_A$ generated by $B_A$ and all $E_i^{(n)}$, $n \in \N$, see \cite{APW}. We
let $H^0_{A, i}$ denote the induction functor from $\mathcal C_A^-$ to
$\mathcal C(P_A(i))$
where this last category consists of all integrable $P_A(i)$-modules.
The right derived functors are denoted $H_{A,i}^j$, $j \geq 0$. Then we have
the following analogue of Lemma 3.6.

\newpage

\begin{lemma}[cf. \cite{APW} Section 4] Let $\lambda \in X$.
\begin{enumerate}
\item [a) ] If $\langle \lambda, \alpha_i^\vee \rangle \geq 0$ then
  $H^j_{A, i} (\lambda) = 0$ for all $j >0$ and $H^0_{A, i} (\lambda)$
  is a free $A$-module whose weights are $\lambda, \lambda - \alpha_i,
  \dots , s_{\alpha_i} (\lambda)$, all occurring with multiplicity
  $1$.
\item [b) ] If $\langle \lambda, \alpha_i^\vee \rangle < -1$ then
  $H^j_{A, i} (\lambda) = 0$ for all $j \not = 1$ and $H^1_{A, i}
  (\lambda)$ is a free $A$-module whose weights are $\lambda +
  \alpha_i, \lambda + 2 \alpha_i, \dots , s_{\alpha_i} \cdot \lambda$,
  all occurring with multiplicity~$1$.
\item [c) ] If $\langle \lambda, \alpha_i^\vee \rangle = r \geq 0$
  then $\Hom_{P_{A}(i)}(H^1_{A, i} (s_{\alpha_i} \cdot \lambda),
  H^0_{A, i} (\lambda)) \simeq A$. Moreover, $H^0_{A, i} (\lambda)$,
  respectively $H^1_{A, i} (s_{\alpha_i} \cdot \lambda)$ has a
  standard $A$-basis $\{v_0, v_1, \dots , v_r\}$, respectively
  $\{v_0', v_1', \dots , v_r'\}$ with $v_j$, respectively $v'_j$,
  having weight $\lambda - j \alpha_i$, $j = 0, 1, \dots, r$. A
  generator $c_i(\lambda)$ of $\Hom_{P_{A, i}}(H^1_{A, i}
  (s_{\alpha_i} \cdot \lambda), H^0_{A, i}(\lambda)) $ is given by
  (with $(d_1, d_2, \dots , d_n)$ being a minimal $n$-tuple in $\N$
  making the Cartan matrix for $R$ symmetric)
$$c_i(\lambda)(v'_j) = \genfrac{[}{]}{0pt}{}{r}{j}_{d_i} v_j, \ j= 0, 1, \dots , r.$$
\end{enumerate}
\end{lemma}
The gaussian binomial coefficients $ \genfrac{[}{]}{0pt}{}{r}{j}_{d_i}$ occurring in c)
are defined like the usual binomial numbers with each integer $m \in \N$ replaced by
\smash{$[m]_{d_i} = \frac{v^{d_im} - v^{-d_im}}{v^{d_i} - v^{-d_i}}$}. If $d_i = 1$ we omit
this subscript.

\subsection{Divisors}
Let $\mathcal D(A)$ denote the divisor group for $A$, i.e., the free $\Z$-module with
basis consisting of all irreducible polynomials in $A$ with leading coefficients equal to $1$.
If $a \in A \setminus \{0\}$ then we write $\di (a) \in \mathcal D(A)$ for the divisor
associated with $a$. The coefficient corresponding to $v-q$ in $ \di (a)$ we shall
denote $\di _q (a)$.

The following formulas are easy exercises, see \cite{T}, Lemma 5.2:
If $\cha k = 0$ we have
\begin{eqnarray}
\di _q([m]) =
\begin{cases}
1 &\text {if } l \text { divides } m,\\
0 &\text { otherwise.}
\end{cases}
\end{eqnarray}

If $\cha k = p>0$ we have
\begin{eqnarray}
\di _q([m]) =
\begin{cases}
p^{\nu_p(m)} &\text {if } l \text { divides } m,\\
0 &\text { otherwise.}
\end{cases}
\end{eqnarray}

Also if $M$ is a finitely generated torsion $A$-module then $M \simeq
\bigoplus_i A/(a_i)$ for some $a_i \in A$ and we write $\di (M) = \sum_i \di (a_i)$. Clearly,
$\di $ is then additive on short exact sequences of finitely generated torsion $A$-modules.
Again $\di _q(M)$ picks out the coefficient in $\di (M)$ corresponding to $v-q$.

If $M$ is a $U_A^0$-module which is a direct sum of its weight spaces $M_\mu$, and if
$M$ is a finitely generated torsion $A$-module then we define
$\di _{U^0} (M) \in \mathcal D(A)[X]$
by
$$ \di_{U^0} (M) = \sum_{\mu \in X} \di (M_\mu)e^\mu.$$

\subsection{Euler type formulas}
We have now reached the point where we can just mimic what we did in Chapter 4. In particular,
for each $\lambda \in X^+$ and for any $V \in \mathcal C_A$ which is a finitely generated torsion $A$-module we define
$$ e_\lambda^U(V) = \sum_{i \geq 0} (-1)^i \di (\Ext^i_{\mathcal C_A}(\Delta_A(\lambda), V)).$$
Note that just as we had finiteness results for $H_A^j$ in Section 6.5 we also have such results for
$\Ext^j_{\mathcal C_A}$ so that this definition makes sense.

Likewise if $M \in \mathcal C^-_A$ is a finitely generated torsion $A$-module then
$$ e_\lambda^B(M) = \sum_{j\geq0} (-1)^j e_\lambda^U(H^j_A(M)).$$

Then the direct analogue of Theorem 4.1 holds with the same proof: for a) we first reduce to the
case (corresponding to (4.2)) where $M$ is determined by the exact sequence
$$ 0 \rightarrow A_\mu \stackrel{a} \rightarrow A_\mu \rightarrow M \rightarrow 0$$
with $\mu \in X$ and $a \in A\setminus \{0\}$. Then we proceed as in Section 4.3. Also the 
proof of b) is a direct translation.

The next step is for a fixed $\mu \in X^+$ to compute $e_\lambda^U(Q(\mu))$ where $Q(\mu)$
is the cokernel of the homomorphism $c_\mu$ discussed in Section 6.7. Just as in Section 4.5 we
factorize $c_\mu = \tilde c_1 \circ \tilde c_2 \circ \cdots \circ \tilde c_N$ (relative to some
reduced decomposition of $w_0$) and then proceed as in Section 4.6.

Lemma 6.1 c) tells us that the cokernel $Q_i(\lambda)$ of $c_i(\lambda)$ has weights
$\lambda - \alpha_i, \lambda - \alpha_i, \dots , s_{\alpha_i}\cdot \lambda$ and that
the weight space $Q_i(\lambda)_{\lambda - j\alpha_i}$ equals
$A/(\genfrac{[}{]}{0pt}{}{r}{j}_{d_i} )$ with $r = \langle \lambda, \alpha_i^\vee \rangle$.

All this leads exactly as in Sections 4.7-8 to the following
\begin{thm}
Let $\lambda , \mu \in X^+$. The cokernel $Q(\mu)$ of the canonical homomorphism $\Delta_A(\mu)
\rightarrow \nabla_A(\mu)$ satisfies
$$e_\lambda^U(Q(\mu)) = -\sum_{\beta \in R^+} \sum_{(x,m) \in V_\beta(\lambda, \mu)} (-1)^{l(x)}\di ([m]_{d_\beta}).$$
\end{thm}

\subsection{Sum formulas for quantized Weyl modules.}

We now deduce sum formulas by proceeding as in Chapter 5. We first
define for $\mu \in X^+$ a filtration of $\Delta^i_A(\mu)$ of $\Delta_A(\mu)$
by setting
$\Delta^i_A(\mu) = c_\mu^{-1} ((v-q)^i \nabla_A(\mu))$.
The quantized Jantzen's filtration is then the descending filtration of
$\Delta_q(\mu)$
defined by setting $\Delta_q^i(\mu)$ equal to the
image of $\Delta^i_A(\mu)$ under the canonical projection
$\Delta_A(\mu) \rightarrow \Delta_A(\mu) \otimes _A k \simeq \Delta_q(\mu)$.

Taking into account the identities (6.4-5) we now get

\begin{thm} Let $\mu \in X^+$.
\begin{enumerate}
\item [a)] Assume $\cha k = 0$. Then
$$
\sum_{i>0} \ch(\Delta_q^i(\mu)) = \sum_{\beta \in R^+}
\sum_{0 < m < \langle \mu + \rho, \beta^\vee \rangle}
 \chi(\mu - ml \beta).$$

\item [b)] Assume $\cha k = p > 0$. Then 
$$
\sum_{i>0} \ch(\Delta_q^i(\mu)) = \sum_{\beta \in R^+} 
\sum_{0 < m < \langle \mu + \rho, \beta^\vee \rangle}
p^{\nu_p(m)} \chi(\mu - ml \beta).
$$
\end{enumerate}
(As in Chapter 5 $\nu_p$ denotes the $p$-adic valuation.)
\end{thm}

\subsection{Sum formulas for quantized tilting modules.}
Tilting modules for $U_q$ are defined in direct analogy with the way it was 
done for $G$. This means that a finite dimensional $U_q$-module $Q$ is tilting if
it has a filtration where the quotients are Weyl modules 
$\Delta_q(\lambda)$  as well as a filtration where the quotients are  
dual Weyl modules $\nabla_q(\lambda)$. Moreover, each such tilting module
$Q$ has a unique lift to a tilting module $\tilde Q$ for $U_{\tilde A}$ where $\tilde A$
denotes the localization $A_{(v-q)}$ of $A$ at the maximal ideal generated by
$v-q$.

For each $\lambda \in X^+$ we set $\bar F_\lambda(Q) = \Hom_{U_q}(\Delta_q(\lambda), Q)$
and $F_\lambda(\tilde Q) = \Hom_{U_{\tilde A}}(\Delta_{\tilde A}(\lambda), \tilde Q)$. Then 
$F_\lambda (\tilde Q) \otimes _{\tilde A} k \simeq \bar F_\lambda (Q)$.

We define a filtration of  $F_\lambda(\tilde Q)$ consisting of the $\tilde A$-submodules
$$F_\lambda(\tilde Q)^j = \{\phi \in F_\lambda(\tilde Q) \mid 
\psi \circ \phi \in (v-q)^j \tilde A c_\lambda 
\text { for  all } \psi \in \Hom_{U_{\tilde A}}(\tilde Q, \nabla_{\tilde A}(\lambda)) \}.$$ 
The image in $\bar F_\lambda (Q)$ of this filtration is then a $k$-space filtration
whose $j$'th term we denote $\bar F_\lambda (Q)^j$.

Using notation analogous to the one in Chapter 5 we now have the sum formulas, cf. 
\cite{Atiltsum}

\begin{thm} Let $Q$ be a tilting module for $U_q$ and let $\lambda \in X^+$. 
\begin{enumerate}
\item [a)] Assume $\cha k = 0$. Then 
$$ \sum_{j > 0} \dim  {\bar F}_{\lambda}({Q})^j = 
- \sum_{\alpha \in R^+}
\sum_{n < 0 \text{ \rm or } nl > \langle \lambda + \rho, \alpha^\vee \rangle}
[{Q} : \chi(\lambda - nl \alpha)].
$$

\item [b)] Assume $\cha k = p > 0$. Then 
$$ \sum_{j > 0} \dim  {\bar F}_{\lambda}({Q})^j = 
- \sum_{\alpha \in R^+}
\sum_{n < 0 \text{ \rm or } nl > \langle \lambda + \rho, \alpha^\vee \rangle}
p^{\nu_p(n)}[{Q} : \chi(\lambda - nl \alpha)].
$$
\end{enumerate}
\end{thm}

\section{Root subsets. Examples.}

In this section we have collected some remarks and examples concerning the sets
$V(\lambda, \mu)$ and $U(\lambda, \mu)$ occurring in Chapter 4. These sets
play important roles in our proof of the sum formulas. Even though they are defined
in a completely elementary way they are somewhat complicated to describe explicitly.
Fixing distinct dominant weights $\lambda$ and $\mu$ we will explore the implications
of the following key condition involved in the definition of these sets.
\begin{align}
\lambda - n \gamma = w \cdot \mu  \text{ for some } n \in \Z, \gamma \in R \text{ and } w \in W.
\end{align}
Note that this condition is symmetric in $\lambda$ and $\mu$. Further, by switching the 
signs of $n$ and $\gamma$ if necessary, we may require $\gamma \in R^+$, but we prefer 
not to do so here. Instead we define, again for distinct $\lambda, \mu \in X^+$,
\begin{align}
S(\lambda, \mu) = \{\gamma \in R^+ \mid  \lambda - n \gamma = w \cdot \mu  
\text{ for some } n \in \Z \text{ and } w \in W\}.
\end{align}

\subsection{Alternative descriptions of $V(\lambda, \mu)$ and $U(\lambda, \mu)$.}
We will show that (7.1) forces $\lambda < \mu$ or $\mu < \lambda$, leading to the
following result.

\begin{prop}
For $\lambda,\mu \in X^+$ the sets $U(\lambda,\mu)$ and $V(\lambda,\mu)$ are empty unless
$\lambda < \mu$. If $\lambda < \mu$, then
\begin{align*}
U(\lambda,\mu) &= \{(\gamma,w,n)\mid \gamma \in R^+ , w \in W, n \in
\Z, \lambda - n \gamma = w \cdot \mu \},
\\
V(\lambda,\mu) &= \{(\beta,x,m)\mid \beta \in R^+ , x \in W, m \in
\Z, \mu - m \beta = x \cdot \lambda \}.
\end{align*}
\end{prop}

\pf
Throughout the proof suppose that $w \cdot \mu = \lambda - n \gamma$ with $\gamma \in R^+$.
So $\langle \mu + \rho , w^{-1}\gamma^\vee \rangle = \langle \lambda + \rho, \gamma^\vee \rangle - 2n$
as in (4.6). Note that $\lambda,\mu \in X^+$ implies $w \cdot \mu \leq \mu$ and
$w^{-1}\cdot \lambda \leq \lambda$ with equalities iff $w$ is the identity.

It will be convenient to prove the claims for $U(\lambda,\mu)$ and $V(\mu, \lambda)$
simultaneously. (Note the interchanged roles of $\lambda$ and $\mu$ in the latter set.)
It suffices to show that one must have $\lambda < \mu$, $\lambda = \mu$ or $\mu < \lambda$
appropriately depending on the value of $n$.

{\it Case I.} If $n < 0$, then $\lambda < \lambda - n \gamma = w \cdot \mu \leq \mu$.

{\it Case II.} If $n = 0$, then $w =$ the identity and $\mu = \lambda$.

{\it Case III.} If $0 < n < \frac{1}{2} \langle \lambda + \rho , \gamma^\vee\rangle$,
then $\langle \mu + \rho , w^{-1}\gamma^\vee \rangle > 0$. So $w^{-1}\gamma \in R^+$
and hence $\mu < \mu + n(w^{-1}\gamma) = w^{-1} \cdot \lambda \leq \lambda$.

$n = \frac{1}{2} \langle \lambda + \rho , \gamma^\vee\rangle$ is impossible, e.g.,
because that would mean $\langle \mu + \rho , w^{-1}\gamma^\vee \rangle = 0$.

To deal with the remaining possibilities, we use
$s_{\gamma}w \cdot \mu = \lambda - (\langle \lambda + \rho , \gamma^\vee \rangle - n)\gamma$
(see Remark 4.7). If
$\frac{1}{2} \langle \lambda + \rho, \gamma^\vee\rangle < n < \langle \lambda + \rho , \gamma^\vee\rangle$,
one reduces to Case III and concludes that $\mu < \lambda$.
Similarly if $n = \langle \lambda + \rho , \gamma^\vee\rangle$, then $\mu = \lambda$
via Case II and if $n > \langle \lambda + \rho , \gamma^\vee\rangle$, then $\lambda < \mu$
via Case I.

\subsection{Explicit determination of the sets $U(\lambda,\mu)$ and $V(\lambda,\mu)$.}
We start by making some easy reductions towards computing these sets.
First, we remark that it is enough to determine the sets $S(\lambda,
\mu)$ defined in (7.2).  By the proof of Proposition~7.1, $S(\lambda,
\mu)$ is nonempty precisely when $U(\lambda,\mu) \cup V(\mu, \lambda)$
is nonempty. If this happens, exactly one of the sets in the union is
nonempty (depending on whether $\lambda < \mu$ or $\mu < \lambda$).
Then, by Remark 4.7, the size of this set is $2|S(\lambda, \mu)|$.
Further, we will see that in all our examples, for each $\gamma \in
S(\lambda, \mu)$, the two associated values of $n$ (and the
corresponding $w \in W$) are easy to determine.

Next, we reduce to the case of irreducible root systems. Note that $S(\lambda, \mu)$ is
described directly in terms of the root system $R$. Clearly, for (7.1) to hold,
$\lambda$ and $\mu$ must differ only in the component of $R$ to which $\gamma$
belongs. In particular the largest possible cardinality of $S(\lambda, \mu)$
for $R$ is the maximum of this cardinality for the irreducible components of
$R$. In the remaining sections we will describe all the different possibilities
that can occur for (simply connected almost simple) groups of classical types
$A, B, C$ and $D$. We start by summarizing part of the findings.

\begin{prop}
When nonempty, the sets $V(\lambda, \mu)$ and $U(\lambda, \mu)$ have
cardinality $2$ for type $A_m$, cardinality $2$ or $4$ for types
$D_m$ ($m > 3$) and $B_2$, and cardinality $2$, $4$ or $6$ for types
$B_m$ and $C_m$ ($m > 2$).
\end{prop}

Looking at (7.1), it makes sense that the sets in question are smaller for
sparser root systems. For type $G_2$ one can check that the cardinality of
these sets is again 0, 2, 4 or 6. We did not work out the types $F_4, E_6, E_7$
and $E_8$.

\subsection{Notation} Let us fix some notation that will be in force
throughout the remaining sections. We will realize the classical root
systems in standard ways (recalled below) in $\R^m$. We will use a fixed
orthonormal basis $\{\epsilon_1, \epsilon_2, \dots , \epsilon_m\}$ for $\R^m$.
Then the set of weights $X$ is a subset of $\R^m$. For any $\lambda \in X$ we set
$$\lambda + \rho = \sum_{i=1}^m \lambda_i \epsilon_i \quad
\text{and } \quad I_\lambda = \{\lambda_1, \dots, \lambda_m\} .$$
Moreover, if also $\mu \in X$ we define the difference sets
$$D_{\lambda \mu} = I_\lambda \setminus I_\mu \quad \text{ and } \quad
D_{\mu \lambda} = I_\mu \setminus I_\lambda.$$
In the following we fix two dominant weights $\lambda$ and $\mu$ and for the classical
types we describe the set $S(\lambda, \mu)$. Since
$I_\lambda$ and $I_\mu$ both have cardinality $m$, we get
$|D_{\lambda \mu}| = |D_{\mu \lambda}|$. Nonemptiness of $S(\lambda, \mu)$
will be characterized by these set differences having cardinality exactly
1 or 2 along with some easy numerical conditions.

\subsection{Type $A$} For $m>1$, we realize $R$ of type $A_{m-1}$ as
the subset $R= \{\epsilon_i - \epsilon_j \mid i \neq j, 1 \leq i, j
\leq m \}$ of $\R^{m}$ with positive roots defined by the condition $i
< j$.  A~vector $\sum_{i=1}^{m} q_i \epsilon_i$ is a weight precisely
when each $q_i - q_j \in \Z$ and $\sum_{i=1}^{m} q_i = 0$ (i.e., when,
for some $t \in \Z$, every $q_i \in \frac{t}{m} + \Z$.) The Weyl group
acts on $\R^m$ by permuting the $\epsilon_i$. So using the notation
from Section 7.3 we see that the $W$-orbit (under the `dot' action) of
any weight $\eta$ consists of those $\eta ' \in X$ for which $I_\eta =
I_{\eta '}$.  Note that $\eta \in X^+$ is equivalent to the condition
$\eta_1 > \eta_2 > \dots > \eta_m$.

Now assume that (7.1) holds for (our fixed and distinct $\lambda, \mu \in X^+$ and) $\gamma = \epsilon_a - \epsilon_b$.
This is equivalent to having
$$I_\mu = (I_\lambda \setminus \{\lambda_a, \lambda_b\}) \cup \{\lambda_a - n, \lambda_b + n\},$$
which implies that $D_{\lambda\mu} = \{\lambda_a, \lambda_b\}$. Thus
a positive $\gamma$ is uniquely determined by $\lambda$ and $\mu$. Morover,
if $D_{\mu\lambda} = \{\mu_c, \mu_d\}$, then one must have
$n = \lambda_a - \mu_c = \mu_d - \lambda_b$ or
$n = \lambda_a - \mu_d = \mu_c - \lambda_b$.

Conversely, we always have $|S(\lambda, \mu)| = 0$ or $1$. The
latter occurs precisely when $|D_{\lambda \mu}| = |D_{\mu \lambda}| = 2$.
In that case, letting $D_{\lambda\mu} = \{\lambda_a, \lambda_b\}$ and
$D_{\mu\lambda} = \{\mu_c, \mu_d\}$, we automatically have
$\lambda_a + \lambda_b = \mu_c + \mu_d$ (since
$\sum_{i=1}^{m} \lambda_i = \sum_{i=1}^{m} \mu_i = 0$).  Then the
corresponding values of $n$ can be read off as above and the
associated permutations $w$ are also easy to describe explicitly.

Suppose $G = GL_m$ instead of $SL_m$. Take $\lambda$ and $\mu$ to
be partitions with at most $m$ parts. (This just gives a different
language to address the question at hand without altering it--the
equivalence is given via translation by a (possibly fractional)
multiple of the $W$-invariant vector $\sum_{i=1}^m \epsilon_i$.)
Then it turns out that $|S(\lambda, \mu)| = 1$ precisely when the
Young diagrams of $\lambda$ and $\mu$ ``differ by connected skew
hooks'', see \cite{Kul1}.

\subsection{Type $B$} For $m>1$, we realize $R$ of type $B_m$ as the subset
$R= \{\pm \epsilon_i, \pm \epsilon_i \pm \epsilon_j \mid i \neq j, 1 \leq i, j \leq m \}$
of $\R^{m}$. Positive roots are those of the form either $\epsilon_i \pm \epsilon_j$
with $i < j$ or $\epsilon_i$.  A vector $\sum_{i=1}^{m} q_i \epsilon_i$
is a weight precisely when each $q_i \in \Z$ or each $q_i \in \frac{1}{2} + \Z$. 
The Weyl group acts on $\R^m$ by permuting the $\epsilon_i$ and 
by changing the signs of $\epsilon_i$. Hence two weights $\eta$ and $\eta'$ belong to 
the same $W$-orbit (under the `dot' action) if and only if the two sets $I_\eta$ and
$I_{\eta '}$ coincide up to signs. Using the notation in Section 7.3,
a weight $\eta$ is dominant if and only if the condition 
$\eta_1 > \eta_2 > \dots > \eta_m > 0$ is satisfied. 

Assume that (7.1) holds for $\lambda$ and $\mu$. To analyze the implications 
of this, we separate into three cases depending on the form of the root $\gamma$ 
in question.

{\it Case 1.} $\gamma = \epsilon_a$. Then (7.1) is equivalent to
having
$$I_\mu = (I_\lambda \setminus \{\lambda_a\}) \cup \{|\lambda_a - n|\},$$
which implies that $D_{\lambda\mu} = \{\lambda_a\}$. Moreover,
if $D_{\mu\lambda} = \{\mu_c\}$, then we have
$|\lambda_a - n| = \mu_c$, i.e., $n = (\lambda_a \pm \mu_c)$.
On the other hand, Case 1 clearly arises whenever
$|D_{\lambda \mu}| = |D_{\mu \lambda}| = 1$ since, for $m>1$, this 
guarantees that $n = (\lambda_a \pm \mu_c) \in \Z$. The associated 
$w$ are easily deduced.

{\it Case 2.} $\gamma = \epsilon_a - \epsilon_b$ with $a \neq b$. (It is convenient
not to assume $\gamma$ to be positive.) Then (7.1) is equivalent to having
$$I_\mu = (I_\lambda \setminus \{\lambda_a, \lambda_b\}) \cup \{|\lambda_a - n|, |\lambda_b + n|\}.$$
There are now two possibilities.

{\it Subcase 2.1.} $|D_{\lambda\mu}| = 1$. By reversing the signs of $\gamma$ and $n$,
we may assume without loss of generality that $D_{\lambda\mu} = \{\lambda_a\}$.
Suppose $D_{\mu\lambda} = \{\mu_c\}$. Then one of the following two options must be
true. 
\begin{enumerate}
\item [(i)] $\lambda_b + n = - \lambda_b$ and $\lambda_a - n = \pm \mu_c$, or
\item [(ii)]  $\lambda_b + n = \pm \mu_c$ and $\lambda_a - n = \pm \lambda_b$.
\end{enumerate}
By adding, in both cases we have $\lambda_b + \lambda_a = \pm \lambda_b \pm \mu_c$.
By positivity constraints and since $\mu_c \neq \lambda_a$ by assumption,
we must have $\lambda_b + \lambda_a = - \lambda_b + \mu_c$, i.e.,
$\lambda_a - \mu_c = -2\lambda_b$.

Conversely, clearly Subcase 2.1 arises exactly when both the
following conditions hold: $|D_{\lambda \mu}| = |D_{\mu \lambda}| = 1$
and (letting $D_{\lambda\mu} = \{\lambda_a\}$ and $D_{\mu\lambda} = \{\mu_c\}$) 
there exists a necessarily unique $b$ such that $\lambda_a - \mu_c = -2\lambda_b$.
The corresponding two values of $n$ are easily found to be $n = -2\lambda_b = \lambda_a - \mu_c$
(leading to (i)) and $n = \mu_c  - \lambda_b = \lambda_a + \lambda_b$
(leading to (ii)). The associated $w$ are easily deduced.

{\it Subcase 2.2.} $|D_{\lambda\mu}| = 2$. Clearly
$D_{\lambda\mu} = \{\lambda_a, \lambda_b\}$. Let
$D_{\mu\lambda} = \{\mu_c, \mu_d\}$ with $c < d$. Then one of the following two
options must be true. Either $\lambda_b + n = \pm \mu_d$ and $\lambda_a - n = \pm \mu_c$,
or $\lambda_b + n = \pm \mu_c$ and $\lambda_a - n = \pm \mu_d$. Adding and using
$\mu_c > \mu_d$, in both cases we have $\lambda_a + \lambda_b = \mu_c \pm \mu_d$.

Conversely, Subcase 2.2 arises exactly when the following conditions
hold: $|D_{\lambda \mu}| = |D_{\mu \lambda}| = 2$ and (setting 
 $D_{\mu, \lambda}= \{\mu_c, \mu_d\}$ with $c<d$) one of the equalities
$\lambda_a + \lambda_b = \mu_c \pm \mu_d$ holds. Then the corresponding two
values of $n$ and the associated $w$ are easily deduced. 

{\it Case 3.} $\gamma = \epsilon_a + \epsilon_b$ with $a \neq b$. The analysis is
very similar to Case 2, so we skip some details.

{\it Subcase 3.1.} $|D_{\lambda\mu}| = 1$. Suppose $D_{\lambda\mu} = \{\lambda_a\}$
and $D_{\mu\lambda} = \{\mu_c\}$. As in Subcase~2.1, we deduce
$\lambda_b - \lambda_a = - \lambda_b \pm \mu_c$, i.e.,
$\lambda_a \pm \mu_c = 2\lambda_b$.

Conversely, clearly Subcase 3.1 arises exactly when both the
following conditions hold. First, $|D_{\lambda \mu}| = |D_{\mu \lambda}| = 1$.
Next, there must exist a necessarily unique $b$ such that $\lambda_a - \mu_c = 2\lambda_b$
and/or a necessarily unique $b'$ such that $\lambda_a + \mu_c = 2\lambda_{b'}$.
If $\lambda_a - \mu_c = 2\lambda_b$, then $n = 2 \lambda_b = \lambda_a - \mu_c$ 
or $n = \lambda_a - \lambda_b = \lambda_b + \mu_c$. If $\lambda_a + \mu_c = 2\lambda_{b'}$, then 
$n = 2 \lambda_{b'} = \lambda_a + \mu_c$ or $n = \lambda_a - \lambda_{b'} = \lambda_{b'} - \mu_c$.

{\it Subcase 3.2.} $|D_{\lambda\mu}| = 2$. Clearly
$D_{\lambda\mu} = \{\lambda_a, \lambda_b\}$. Suppose $a < b$. Let
$D_{\mu\lambda} = \{\mu_c, \mu_d\}$ with $c < d$. As before, we deduce
$\lambda_a - \lambda_b = \mu_c \pm \mu_d$.

Conversely, Subcase 3.2 arises exactly when the following conditions
hold. $|D_{\lambda \mu}| = |D_{\mu \lambda}| = 2$ and one of the equalities
$\lambda_a - \lambda_b = \mu_c \pm \mu_d$ holds. Again the corresponding two
values of $n$ and the associated $w$ are easily deduced. 

Using the above cases, we sketch a procedure to calculate the sets 
$S(\lambda,\mu)$ (and hence the sets $U(\lambda,\mu)$ and $V(\mu,\lambda)$). 
Clearly $S(\lambda,\mu)$ is empty unless $|D_{\lambda \mu}| = |D_{\mu \lambda}| = 1 \text{ or } 2.$

Suppose $D_{\lambda\mu} = \{\lambda_a, \lambda_b\}$ with $a < b$ and
$D_{\mu\lambda} = \{\mu_c, \mu_d\}$ with $c < d$. Then at most one of cases 2.2
and 3.2 can occur, since $\lambda_a + \lambda_b = \mu_c \pm \mu_d$ 
and $\lambda_a - \lambda_b = \mu_c \pm \mu_d$ cannot be true simultaneously. 
So $|S(\lambda,\mu)| = 0$ or 1.

Suppose $D_{\lambda\mu} = \{\lambda_a\}$ and $D_{\mu\lambda} = \{\mu_c\}$.
Then Case 1 will always occur and cases 2.1 and 3.1 will occur depending on the
existence of $b$ satisfying one of the three conditions $\lambda_a - \mu_c = -2\lambda_b$, 
$\lambda_a - \mu_c = 2\lambda_b$ and $\lambda_a + \mu_c = 2\lambda_b$. 
Clearly at most two of these can be satisfied (since at most one of the
first two can be true), each by a unique $b$. Note also that if the rank
$m = 2$, then at most one of the three conditions can hold. So $|S(\lambda,\mu)| = 3$
(provided $m > 2$) or 2 or 1.

An example with $|S(\lambda,\mu)| = 3$ for type $B_3$ is given by 
$\lambda + \rho = 5 \epsilon_1 + 3 \epsilon_2 + 2 \epsilon_3$ 
and $\mu + \rho = 3 \epsilon_1 + 2 \epsilon_2 + \epsilon_3$.

\subsection{Type $C$} The analysis can be lifted almost verbatim from
that for type $B$, so we indicate only the changes that need to be
made there.  In $R$ we replace $\pm \epsilon_i$ by~$\pm 2\epsilon_i$,
with $2\epsilon_i \in R^+$.  For weights we require each $q_i \in \Z$.
Again we make three cases.  Only Case 1 needs any change. Here we take
$\gamma = 2\epsilon_a$. Then (7.1) leads to $n = \frac{1}{2}(\lambda_a
\pm \mu_c)$. Conversely, this case arises exactly when $|D_{\lambda
  \mu}| = |D_{\mu \lambda}| = 1$ and $(\lambda_a \pm \mu_c)$ is even.
Except for the inclusion of the evenness condition, the procedure to
calculate $S(\lambda,\mu)$ stays unchanged. (In particular
$|D_{\lambda \mu}| = |D_{\mu \lambda}| = 1$ no longer guarantees
$|S(\lambda,\mu)| \geq 1$.)

\subsection{Type D} The analysis is again similar to that for type $B$, 
so we indicate only the changes. Here they are more significant. Now 
$R= \{\pm \epsilon_i \pm \epsilon_j \mid i \neq j, 1 \leq i, j \leq m \}$
with positive roots those of the form $\epsilon_i \pm \epsilon_j$ with 
$i < j$.  The weights stay the same, i.e., $\sum_{i=1}^{m} q_i \epsilon_i$ 
with each $q_i \in \Z$ or each $q_i \in \frac{1}{2} + \Z$. The Weyl group acts 
by permuting the $\epsilon_i$ and by changing the signs of an even number 
of $\epsilon_i$. This means that two weights $\eta$ and $\eta'$ belong to the same 
$W$-orbit (under the `dot' action) if and only if the two sets $I_\eta$ and $I_{\eta'}$
coincide up to an even number of signs. In this case we shall therefore find it
convenient to  work with the set $I'_\eta = \{|\eta_1|, |\eta_2|, \dots , |\eta_m|\}$
instead of $I_\eta$, and we replace $D$ by $D'$ for the corresponding difference sets.
Using otherwise the notation in Section 7.3, $\eta \in X^+$
is equivalent to the condition 
$\eta_1 > \eta_2 > \dots > \eta_{m-1} > |\eta_m|$. 
Note that $\eta_m$ may be $0$ or negative. If $\eta \in X^+$ with $\eta_m \geq 0$ then $I'_\eta = I_\eta$.

Assuming (7.1), for $\gamma = \epsilon_a \pm \epsilon_b$,
we get
$$I'_\mu = (I'_\lambda \setminus \{|\lambda_a|, |\lambda_b|\}) \cup
\{|\lambda_a - n|, |\lambda_b \mp n|\}.$$
Now just as for type $B$, we get the following consequences.

If $\gamma = \epsilon_a + \epsilon_b$ (respectively, $\epsilon_a -
\epsilon_b$)
and $D'_{\lambda\mu} = \{|\lambda_a|, |\lambda_b|\}$, then letting
$a<b$ and $D'_{\mu\lambda} = \{\mu_c, |\mu_d|\}$ with $c<d$, we get
$\lambda_a - \lambda_b = \mu_c \pm \mu_d$ (respectively,
$\lambda_a + \lambda_b = \mu_c \pm \mu_d$).

If $\gamma = \epsilon_a + \epsilon_b$ (respectively, $\epsilon_a -
\epsilon_b$) and $D'_{\lambda\mu} = \{|\lambda_a|\}$, then letting
$D'_{\mu\lambda} = \{|\mu_c|\}$, we get $\lambda_a \pm  \mu_c = 2 \lambda_b$ 
(respectively, $\lambda_a - \mu_c = - 2 \lambda_b$; unlike for type $B$,
here it requires some work to rule out $\lambda_a + \mu_c = -2\lambda_b$. One
sees that the latter equality only arises when $\mu_c = 0$, in which case one may as well
use $\lambda_a - \mu_c = -2\lambda_b$).

Conversely, we now describe exactly when (7.1) holds for a given $\gamma$. 
We make the convention that $\sign(0) = 0$.
For a weight $\eta$, define $\sign(\eta) = \Pi_{i=1}^m \sign(\eta_i)$.
Note that for all $\eta$ in a $W$-orbit (under the `dot' action), $\sign(\eta)$ 
remains the same. Clearly, making cases as for Type B, the validity of (7.1) 
in each case is characterized by the respective numerical constraints along
with the requirement $\sign(\mu) = \sign(\lambda - n \gamma)$. We make this 
explicit below. Since $\mu \in X^+$, $\sign(\mu) = \sign(\mu_m)$ (and so the 
sign condition is vacuous if $\mu_m =0$). To calculate $\sign(\lambda - n \gamma)$, 
we have used the values of $n$ obtained in each case.

Suppose $|D'_{\lambda \mu}| = |D'_{\mu \lambda}| = 2$. Let
$D'_{\lambda\mu} = \{\lambda_a, |\lambda_b|\}$ with $a<b$
and $D'_{\mu\lambda} = \{\mu_c, |\mu_d|\}$ with $c<d$.
Then (7.1) holds for $\gamma = \epsilon_a + \epsilon_b$
iff $\lambda_a - \lambda_b = \mu_c \pm \mu_d$ and
$$ \sign(\mu_m) = \begin{cases}  \sign(-\lambda_a + \lambda_b + \mu_c) \sign (\lambda_m) &\text { if } b<m;\\
\sign(-\lambda_a + \lambda_m + \mu_c) &\text { if } b = m.\end{cases}$$
Similarly, (7.1) holds for $\gamma = \epsilon_a - \epsilon_b$ iff
$\lambda_a + \lambda_b = \mu_c \pm \mu_d$ and
$$ \sign(\mu_m) = \begin{cases}  \sign(\lambda_a + \lambda_b - \mu_c) \sign (\lambda_m) &\text { if } b<m;\\
\sign(\lambda_a + \lambda_m - \mu_c) &\text { if } b = m.\end{cases}$$
Note that in the above cases, respectively, $-\lambda_a + \lambda_b + \mu_c = \mp \mu_d$
and $\lambda_a + \lambda_b - \mu_c =\nobreak  \pm \mu_d$.

Suppose $|D'_{\lambda \mu}| = |D'_{\mu \lambda}| = 1$. Let
$D'_{\lambda\mu} = \{|\lambda_a|\}$  and $D'_{\mu\lambda} = \{|\mu_c|\}$.
Then (7.1) holds for $\gamma = \epsilon_a +\epsilon_b$ iff
$\lambda_a \pm \mu_c = 2\lambda_b$ and
$$ \sign(\mu_m) = \begin{cases}  - \sign(\lambda_a - 2\lambda_b) \sign (\lambda_m) &\text { if } a<m;\\
1 &\text { if } a = m.\end{cases}$$
Likewise (7.1) holds for $\gamma = \epsilon_a -\epsilon_b$ iff
$\lambda_a - \mu_c = - 2\lambda_b$ and
$$ \sign(\mu_m) = \begin{cases}  - \sign(\lambda_a + 2 \lambda_b) \sign (\lambda_m) &\text { if } a<m;\\
-1  &\text { if } a = m.\end{cases}$$
Again note that in these cases we have, respectively, 
$\lambda_a - 2 \lambda_b = \mp \mu_c$ and $\lambda_a + 2\lambda_b =
\nobreak \mu_c$. 

By easy extensions of the arguments for type B, one easily deduces
the following. If $|D'_{\lambda \mu}| = |D'_{\mu \lambda}| = 2$,
then $|S(\lambda, \mu)| \leq 1$. If $|D'_{\lambda \mu}| = |D'_{\mu \lambda}| = 1$,
then $|S(\lambda, \mu)| \leq 2$. In fact $|S(\lambda, \mu)| = 2$
can occur only when $\mu_m = \lambda_m = 0$. This can be seen using the sign
constraints above. In particular, for types $D_2$ and $D_3$, $|S(\lambda, \mu)| \leq 1$
(as we already know from the result for type $A$). For type $D_4$, the example
given for type $B_3$ provides an instance where $|S(\lambda, \mu)| = 2$.


\begin{thebibliography}{99}

\bibitem{Asum}
Henning Haahr Andersen, {\it Filtrations of cohomology modules for Chevalley
groups},
{ Ann.\ Scient.\ \'Ec.\ Norm.\ Sup.\ }(4) {\bf 16} (1983),
495--528.

\bibitem{APW}
Henning Haahr Andersen, Patrick Polo and Wen Kexin, {\it Representations of
quantum algebras}, Invent. Math. {\bf 104} (1991), 1--59.

\bibitem{AW}  
Henning Haahr Andersen and Wen Kexin, {\it Representations of quantum algebras. The mixed case},
J. reine angew. Math. {\bf 427} (1992), 35--50

\bibitem{Atilt}Henning Haahr Andersen, {\it Filtrations and tilting modules},
{ Ann.\ Scient.\ \'Ec.\ Norm.\ Sup.\ }(4) {\bf 30} (1997),
353--366.

\bibitem{Atiltsum}Henning Haahr Andersen, {\it A sum formula for tilting 
filtrations}. Commutative algebra, homological algebra and representation theory (Catania/Genoa/Rome, 1998).  
J. Pure Appl. Algebra  152  (2000),  no. 1--3, 17--40

\bibitem{Bott}
Raoul Bott, {\it Homogeneous vector bundles}, 
Ann. of Math. (2) {\bf 44} (1957), 203--248.

\bibitem{Dem}
Michel Demazure, {\it  A very simple proof of Bott's theorem},
Invent. math. {\bf 33} (1976), 271--272.

\bibitem{Do}
Stephen Donkin, {\it Rational representations of algebraic groups. 
Tensor products and filtration.} 
Lecture Notes in Mathematics, 1140. Springer-Verlag, Berlin, 1985. vii+254 pp.

\bibitem{RH}
Steen Ryom-Hansen, {\it A $q$-analogue of Kempf's vanishing theorem.}  
Mosc. Math. J.  3  (2003),  no. 1, 173--187, 260.


\bibitem{JCJ}
Jens Carsten Jantzen, {\it Darstellungen halbeinfacher Gruppen und 
kontravariante Formen},  J. Reine Angew. Math. {\bf 290}  (1977), 117--141.

\bibitem{MHG}
Jens Carsten Jantzen, {\it Moduln mit einem h\"ochsten Gewicht}, Lecture Notes in 
Mathematics, {\bf 750}. Springer-Verlag, Berlin Heidelberg New York, 1979.

\bibitem{RAG}
Jens Carsten Jantzen, {\it Representations of algebraic groups}. Second edition,
Mathematical Surveys and Monographs, {\bf 107}. American Mathematical Society, Providence, 
RI, 2003.

\bibitem{Ke}
George R. Kempf, {\it Linear systems on homogeneous spaces}, Ann. of Math. (2) {\bf 103}
(1976), no. 3, 557--591.

\bibitem{Kul1}
Upendra Kulkarni, {\it A homological interpretation of Jantzen's sum formula},
Transformation Groups {\bf 11} (2006), no. 3, 517--538.

\bibitem{Kul2}
Upendra Kulkarni, {\it On Jantzen's and Andersen's sum formulas for algebraic 
groups}, preprint.

\bibitem{L}
George Lusztig, {\em Quantum groups at roots of $1$}, Geom. Ded. {\bf 35} (1993), 89--114. 


\bibitem{T}
Lars Thams, {\em Two classical results in the quantum mixed case}, 
J. reine angew. Math. {\bf 436} (1993), 129--153.

\end{thebibliography}
\end{document}